\documentstyle[12pt,newlfont]{amsart}

\newcommand{\nc}{\newcommand}

\nc{\one}{\mbox{\bf 1}}
\nc{\invtensor}{\underset{\leftarrow}{\otimes}}
\nc{\ad}{\operatorname{\mbox{\em ad}}}
\nc{\rk}{\operatorname{rank}}
\nc{\corank}{\operatorname{corank}}
\nc{\Sym}{\operatorname{Sym}}
\nc{\sym}{\operatorname{\mbox{\em sym}}}
\nc{\id}{\operatorname{id}}
\nc{\htt}{\operatorname{ht}}
\nc{\Ker}{\operatorname{Ker}}
\nc{\im}{\operatorname{Im}}
\nc{\osp}{\operatorname{osp}}
\nc{\ssp}{\operatorname{sp}}
\nc{\sgn}{\operatorname{sgn}}
\nc{\F}{\operatorname{F}}

\nc{\out}{\operatorname{out}}
\nc{\supp}{\operatorname{supp}}
\nc{\Hom}{\operatorname{Hom}}
\nc{\End}{\operatorname{End}}

\nc{\Ann}{\operatorname{Ann}}
\nc{\Ind}{\operatorname{Ind}}

\nc{\wt}{\operatorname{wt}}
\nc{\ch}{\operatorname{ch}}
\nc{\Stab}{\operatorname{Stab}}
\nc{\Sch}{{\cal S}\mbox{\em ch}}

\nc{\Spec}{\operatorname{Spec}}
\nc{\Prim}{\operatorname{Prim}}
\nc{\Aut}{\operatorname{Aut}}
\nc{\Fract}{\operatorname{Fract}}
\nc{\gr}{\operatorname{\mbox{\em gr}}}

\nc{\Dglie}{\operatorname{{\cal D}glie}}
\nc{\Dgcoalg}{\operatorname{{\cal D}gcoalg}}
\nc{\Homcoalg}{\operatorname{{\cal H}omcoalg}}
\nc{\Homlie}{\operatorname{{\cal H}omlie}}
\nc{\Dgmod}{\operatorname{{\cal D}gmod}}
\nc{\Dgcomod}{\operatorname{{\cal D}gcomod}}
\nc{\Card}{\operatorname{Card}}
\nc{\pa}{\partial}
\nc{\co}{\cal O}

\nc{\Ug}{{\cal U}({\frak g})}
\nc{\Zg}{{\cal Z}({\frak g})}
\nc{\Uh}{{\cal U}({\frak h})}
\nc{\Ut}{{\cal U}({\frak t})}
\nc{\fg}{\frak g}
\nc{\cS}{\cal S}
\nc{\cU}{\cal U}
\nc{\cA}{\cal A}
\nc{\cK}{\cal K}
\nc{\cH}{\cal H}
\nc{\cL}{\cal L}
\nc{\cF}{\cal F}
\nc{\fm}{\frak m}
\nc{\fn}{\frak n}
\nc{\fh}{\frak h}
\nc{\ft}{\frak t}
\nc{\fp}{\frak p}
\nc{\fI}{\frak I}
\nc{\wdM}{\widetilde{M}}
\nc{\wdV}{\widetilde{V}}
\nc{\Sh}{{\cal S}({\frak h})}
\nc{\dirlim}{\underset{\rightarrow}{\lim}\,} 
\nc{\nen}{\newenvironment}
\nc{\ol}{\overline}
\nc{\ul}{\underline}
\nc{\ra}{\rightarrow}
\nc{\lra}{\longrightarrow}
\nc{\Lra}{\Longrightarrow}
\nc{\Lla}{\Longleftarrow}
\nc{\Llra}{\Longleftrightarrow}
\nc{\hra}{\hookrightarrow}
\nc{\iso}{\overset{\sim}{\lra}}
\nc{\ssubset}{\underset{\not=}{\subset}}

\nc{\Thm}[1]{Theorem~\ref{#1}}
\nc{\Prop}[1]{Proposition~\ref{#1}}
\nc{\Lem}[1]{Lemma~\ref{#1}}
\nc{\Cor}[1]{Corollary~\ref{#1}}
\nc{\Conj}[1]{Conjecture~\ref{#1}}
\nc{\Claim}[1]{Claim~\ref{#1}}
\nc{\Defn}[1]{Definition~\ref{#1}}
\nc{\Exa}[1]{Example~\ref{#1}}
\nc{\Rem}[1]{Remark~\ref{#1}}
\nc{\Note}[1]{Note~\ref{#1}}

\nen{thm}[1]{\label{#1}{\bf Theorem.\ } \em}{}
\nen{prop}[1]{\label{#1}{\bf Proposition.\ } \em}{}
\nen{lem}[1]{\label{#1}{\bf Lemma.\ } \em}{}
\nen{cor}[1]{\label{#1}{\bf Corollary.\ } \em}{}
\nen{conj}[1]{\label{#1}{\bf Conjecture.\ } \em}{}
\nen{claim}[1]{\label{#1}{\bf Claim.\ } \em}{}
 
\nen{defn}[1]{\label{#1}{\bf Definition.\ } }{}
\nen{exa}[1]{\label{#1}{\bf Example.\ } }{}

\nen{rem}[1]{\label{#1}{\em Remark.\ } }{}
\nen{note}[1]{\label{#1}{\em Note.\ } }{}
\nen{exer}[1]{\label{#1}{\em Exercise.\ } }{}
\nen{sket}[1]{\label{#1}{\em Sketch of proof.\ } }{}

\begin{document}

%top matter
\title[]{The annihilation theorem for the 
completely reducible Lie superalgebras }

\author[]{Maria Gorelik, Emmanuel Lanzmann}
\address{ Laboratoire de Math\'ematiques Fondamentales, Universit\'e
Pierre et Marie Curie, 4 Place Jussieu, Paris 75252 Cedex 05, France;\\
Dept. of Theoretical Mathematics, The Weizmann Institute of Science, 
Rehovot 76100, Israel,
{\tt email: gorelik@@math.jussieu.fr,\ lanzmann@@wisdom.weizmann.ac.il} 
}

\thanks{
The first author was partially supported by the Hirsch and Braine Raskin
Foundation, Minerva grant  8337 and Chateubriand fellowship}
\thanks{
The second author was partially supported by TMR Grant No. FMRX-CT97-0100 and
Minerva grant 8337}

\begin{abstract}
A well known theorem of Duflo claims that the annihilator of a Verma module 
in the enveloping algebra of
a complex semisimple Lie algebra is generated by its intersection with the 
centre. For a Lie superalgebra this result fails to be true. 
For instance, in the case of the orthosymplectic Lie superalgebra
$\osp(1,2)$, Pinczon gave in~\cite{pin} an example of a Verma module
whose annihilator is not generated by its intersection with the centre of
universal enveloping algebra. More generally,  Musson 
produced  in~\cite{mu} a family of such "singular" Verma modules for
$\osp(1,2l)$ cases.

In this article we give a necessary and sufficient condition on the 
highest weight of a $\osp(1,2l)$-Verma module
for its annihilator to be generated by its intersection with the centre. This 
answers a  question of Musson. 

The classical proof of the Duflo theorem is based on a deep result of 
Kostant which uses
some delicate algebraic geometry reasonings. Unfortunately these arguments 
can not be 
reproduced in the quantum and super cases.
This obstruction
forced Joseph and Letzter, in their work on the quantum case 
(see~\cite{jl}), to find 
an alternative approach to the Duflo theorem. 
Following their ideas, we
compute the factorization of the  Parthasarathy--Ranga-Rao--Varadarajan (PRV)
 determinants.  Comparing it with the factorization
of Shapovalov determinants we find, unlike to the classical and
quantum cases, that the PRV determinant contains some extrafactors. 
The set of 
zeroes of these extrafactors is precisely the set of highest weights of Verma
modules  whose annihilators are not generated by their intersection with the 
centre. 
We also find an analogue of Hesselink formula (see~\cite{he}) giving the 
multiplicity of every simple finite dimensional module in the graded 
component of the harmonic space in the symmetric algebra.

\end{abstract}
\maketitle

\section{Introduction}
The aim of this article is to study the generalization of the annihilator 
theorem of Duflo to 
a certain class of Lie superalgebras, namely to the class of orthosymplectic
 Lie superalgebras $\osp(1,2l)$, $l\geq 1$. The main reason for considering 
this class of Lie superalgebras is the complete reducibility of finite 
dimensional modules. 
Both in classical and quantum cases, the complete reducibility
appears in various  steps of  the study of annihilators of 
Verma modules. A Lie superalgebra whose finite dimensional modules are
completely reducible is called completely reducible. 
According to the theorem of    
Djokovi\'c and Hochschild (see~\cite{sch}, p. 239), any finite dimensional 
completely reducible Lie superalgebra is a direct sum of 
semisimple Lie algebras
and algebras $\osp(1,2l)$ ($l\geq 1$).
Another useful  property of $\osp (1,2l)$ is that its enveloping
algebra is a domain (see~\cite{al}).

Let ${\frak g}={\frak g}_0\oplus {\frak g}_1$ be a Lie superalgebra 
$\osp(1,2l)$, and ${\frak h}$ be a Cartan
subalgebra of the Lie algebra ${\frak g}_0\simeq \ssp(2l)$. Let 
$\Delta=\Delta_0\cup \Delta_1$ be the set of roots of ${\frak g}$
with respect to ${\frak h}$, where
$\Delta_0$ (resp., $\Delta_1$) is the set of even (resp., odd) roots.
Let $\pi$ be a set of simple roots of $\Delta$ and $P^+(\pi)$ be the
set of integral dominant weights.
Denote by $\Delta^+$ (resp.,  $\Delta_0^+,\ \Delta_1^+$)
the set of positive roots (resp., positive even, positive odd).
Set $\overline{\Delta_0^+}\!=\!\Delta_0^+\setminus2\Delta_1^+,\ 
\rho\!=\!{1\over 2}\mathop{\sum}_{\alpha\in \Delta_0^+}\alpha-
{1\over 2}\mathop{\sum}_{\alpha\in \Delta_1^+}\alpha$. Denote
by ${\cal U}({\frak g})$ the universal enveloping algebra
of ${\frak g}$, by ${\cal Z}({\frak g})$ its centre. Let
${\cal F}$ be the canonical filtration of ${\cal U}({\frak g})$.

Consider the case $l=1$, i.e. ${\frak g}=\osp(1,2)$, which has been treated by
Pinczon. 
In this case, any ${\frak g}$-Verma module ${\widetilde M}$,
viewed as a ${\frak g}_0$-module, is the direct sum of two
${\frak g}_0$-Verma modules $M_0$ and $M_1$. Let $C_0$ be 
a Casimir element for ${\frak g}_0$. 
Then $C_0$ acts by
scalars $c_i$ on $M_i$ (i=0,1). In general, $c_1\not=c_0$, 
and in this case, Pinczon proved that 
$\Ann {\widetilde M}={\cal U}({\frak g})\Ann_{{\cal Z}({\frak g})} 
{\widetilde M}$. However when the highest weight of $\widetilde{M}$ equals 
$-\rho$, one has $c_1=c_0$ and so $C_0-c_0$ belongs to the annihilator of 
$\widetilde{M}$. It is easy to check that $(C_0-c_0)\not\in {\cal U}({\frak g})
\Ann_{{\cal Z}({\frak g})}\widetilde{M}$. 
Consequently, the annilihilator theorem 
doesn't hold in full generality.

Return to the general case ${\frak g}=\osp(1,2l)$. For any 
$\lambda\in {\frak h}^*$, let $\widetilde{M}(\lambda)$ be the Verma 
module of highest weight $\lambda$.
Musson showed (see~\cite{mu},5.6) 
that if $(\lambda+\rho,\beta)=0$ for some $\beta\in\Delta_1^+$ and
$(\lambda+\rho,\alpha)\not=0$ for all $\alpha\in \overline{\Delta_0^+}$
then $\Ann_{{\cal U}({\frak g})}\widetilde{M}(\lambda)$ strictly contains
 ${\cal U}({\frak g})\Ann_{{\cal Z}({\frak g})}\widetilde{M}(\lambda)$
and so the annihilator theorem does not hold.
At the end of the article~\cite{mu}, Musson
asked for which $\lambda\in {\frak h}^*$ the annihilator of 
$\widetilde{M}(\lambda)$ 
is generated by its intersection with the centre. 

In this article we give an answer to this question by proving the 
following theorem:

\begin{thm}{} For any $\lambda\in {\frak h}^*$ one has  
\begin{equation}
\label{thmintr} 
\Ann_{{\cal U}({\frak g})}\widetilde{M}(\lambda)=
{\cal U}({\frak g})\Ann_{{\cal Z}({\frak g})}{\widetilde M}(\lambda)\  \ 
\Longleftrightarrow \ \ \forall\alpha\in \Delta_1^+:\ 
(\lambda+\rho,\alpha)\not=0.
\end{equation}
\end{thm}

Our strategy to prove this theorem is based on 
 a recent alternative proof (see~\cite{jfr}) of 
the classical annihilator theorem which  follows essentially the 
quantum case treated by Joseph and Letzter--- see~\cite{jl}. 
This proof  has the following main steps.

1) The first ingredient is the separation theorem established 
by Musson in~\cite{mu}, 1.4. Using complete reducibility, 
he proved the existence of an $\ad{\frak g}$-submodule 
${\cal H}$ of 
${\cal U}({\frak g})$ such that the multiplication map induces an isomorphism
${\cal U}({\frak g})\simeq {\cal Z}({\frak g})\otimes {\cal H}$ of
$\ad{\frak g}$-modules. Moreover, the multiplicity of any finite dimensional 
module in ${\cal H}$ is equal to the dimension of its zero weight space. 
Since the centre acts on any Verma module $\widetilde{M}(\mu)$ 
as a scalar, the separation theorem implies that 
$\Ann_{{\cal U}({\frak g})}\widetilde{M}(\mu)=
{\cal U}({\frak g})\Ann_{{\cal Z}({\frak g})}\widetilde{M}(\mu)$ iff
$\Ann_{{\cal H}}\widetilde{M}(\mu)=0$.

2) In Section~\ref{hesselink}  we obtain a
formula giving the multiplicity of every simple finite dimensional
module $\widetilde{V}(\lambda)$, 
$\lambda\in P^+(\pi)$, in the graded components $H^n$ where 
$H=\gr_{\cal F}{\cal H}$.
This is an analogue of the classical Hesselink formula (see~\cite{he}).
However, in our case the formula looks rather suprising. 
Namely, for every $\lambda\in P^+(\pi)$, the multiplicity of every
$\widetilde{V}(\lambda)$ in ${\cal H}$ is given by the formula
\begin{equation}
[H^n:\widetilde {V}(\lambda)]=\sum_{w\in W} (-1)^{l(w)} P_n(w.\lambda),
\label{superHss}
\end{equation} 
where $P_n(\mu)$ are integers defined by generating function
\begin{equation}
\prod_{\alpha\in \Delta^+_0}(1-q e^{\alpha})^{-1}
(1+q^{2l}e^{\beta_1})\prod_{\beta\in {\Delta^+_1\setminus\{\beta_1\}}}
(1+e^{\beta})=
\sum_{r=0}^{\infty}\sum_{\nu\in {\Bbb N}\pi}P_r(\nu) e^{\nu} q^r
\label{superp}
\end{equation}
where $\beta_1$ is the highest weight of the ${\frak g}_0$-module 
${\frak g}_1$.

3)  Any Verma module $\widetilde{M}(\mu)$ has a simple 
Verma submodule $\widetilde{M}(\mu')$ and $\Ann_{{\cal Z}({\frak g})}
\widetilde{M}(\mu)
=\Ann_{{\cal Z}({\frak g})}\widetilde{M}(\mu')$. Moreover, 
$(\mu+\rho,\alpha)=0$ for some $\alpha\in \Delta_1^+$ iff 
$(\mu'+\rho,\alpha')=0$ for some $\alpha'\in \Delta_1^+$. Thus, in order to
prove the implication "$\Leftarrow$" of~(\ref{thmintr}), it is sufficient
to check it for simple Verma modules. 

For the other implication, we show (see~\ref{otherim}) that there
exists an $\ad {\frak g}$ submodule $V$ of ${\cal U}({\frak g})$, which 
lies in the annihilator of any simple Verma module 
$\widetilde{M}(\mu)$ such that $(\mu+\rho,\alpha)=0$ for some $\alpha\in 
\Delta_1^+$.
Then,~\Lem{zarcl} implies that $V$  lies indeed in the annihilator 
of any Verma module $\widetilde{M}(\mu)$ such that $(\mu+\rho,\alpha)=0$ 
for some $\alpha\in \Delta_1^+$.
So in a sense, both implications of the equivalence~(\ref{thmintr}) are reduced
to the case of a simple Verma module.

4) In Section~\ref{shapo} we present a criterion of simplicity for a 
Verma module
$\widetilde{M}(\mu)$ given by  Kac (see~\cite{k3} and also~\cite{mu2} 2.4). 
As in the classical case, this criterion is related to the so-called  
Shapovalov determinants $\det { S}_{\nu}\in 
{\cal U}({\frak h})$, $\nu\in {\Bbb N}\pi$. These determinants are constructed 
in such a way that $\det {S}_{\nu}(\mu)=0$
means that the Verma module $\widetilde{M}(\mu)$ contains a proper submodule 
with a non trivial element of weight $\mu-\nu$. Thus, $\widetilde{M}(\mu)$ 
is simple iff $\det S_{\nu}(\mu)\not=0$ for all $\nu\in{\Bbb N}\pi$. 
Kac proved
that all  factors of Shapovalov determinants are linear and he gave an explicit
 formula for the 
factorization of these determinants (see~\ref{shapovalov}).

5) The separation theorem allows one to define 
the Parthasarathy--Ranga-Rao--Varadarajan (PRV) determinants as in the 
classical case (see~\ref{prvcons} for details). 
For any $\lambda\in {\frak h}^*$, denote by $\widetilde{V}(\lambda)$
the unique simple module of the highest weight 
$\lambda$.  
The PRV determinant $\det{PRV}^{\lambda}\in {\cal U}({\frak h})$, 
$\lambda\in P^+(\pi)$,  has the property: 
$$\forall \mu\in {\frak h}^*,\ 
\det {PRV}^{\lambda}(\mu)=0\Longleftrightarrow 
\exists \widetilde{V}(\lambda)\subset {\cal H},
\widetilde{V}(\lambda)\subset \Ann_{\cal H}\widetilde{V}(\mu).$$
Thus $\Ann_{{\cal H}}\widetilde{V}(\mu)=0$ iff
$\det{PRV}^{\lambda}(\mu)\not=0$ for all $\lambda\in P^+(\pi).$
Therefore we have to prove that for any simple Verma module
$\widetilde{M}(\mu)=\widetilde{V}(\mu)$ one has the following equivalence
\begin{equation}
\label{prvshap}
\exists \lambda\in P^+(\pi):\ \det{PRV}^{\lambda}(\mu)=0\ \ \Longleftrightarrow
\ \ \exists\alpha\in \Delta_1^+:\ (\mu+\rho,\alpha)=0.
\end{equation}
In Section~\ref{prv} we describe the zeroes of the determinants 
$\det{PRV}^{\lambda}$. 

Both in the classical semisimple and quantum cases
the union of the zeroes of $\det{PRV}^{\lambda}$ ($\lambda\in P^+(\pi)$)
coincides with the union of zeroes of the Shapovalov determinants
$\det { S}_{\nu} (\nu\in {\Bbb N}\pi)$.
This implies the annihilation theorem.
 
6) In~\Thm{thmprv} we give a linear factorization
of the determinants $\det{PRV}^{\lambda}:\ \lambda\in P^+(\pi)$.
There are factors of two types.
The factors of the first type, called "standard factors", coincide
with factors of Shapovalov determinants. They
can be obtained as reviewed in ~\cite{jfr}. We
briefly summarize this procedure in~\ref{stfactor}---~\ref{pfthmLL}.
These factors have form
\begin{equation}
\label{stnd}
\begin{array}{lcc}\varphi(\alpha)+(\rho,\alpha)-m(\alpha,\alpha)/2, &\alpha\in
 \overline{\Delta_0^+},& m\in {\Bbb N}, m\geq 1\\
 \varphi(\alpha)+(\rho,\alpha)-(2m+1)(\alpha,\alpha)/2, &\alpha\in\Delta_1^+, &
m\in {\Bbb N}
\end{array}
\end{equation}
where the element $\varphi(\alpha)$ of ${\frak h}$ is defined by the formula 
$\varphi(\alpha)(\mu)=(\alpha,\mu)$.

A delicate point is to get the factors of the second type, called 
"exotic factors". They have form
\begin{equation}
\label{intrex}
\varphi(\alpha)+(\rho,\alpha),\ \alpha\in \Delta_1^+.
\end{equation}
The Hesselink formula discussed in 2) ensures that  factors (\ref{stnd}), 
(\ref{intrex}) indeed exhaust 
the set of factors of  PRV determinants. 

The exotics factors 
correspond precisely to the equivalence~(\ref{prvshap}). 
We conclude that the annihilator of a simple Verma 
module 
$\widetilde{M}(\mu)=\widetilde{V}(\mu)$ is generated by its intersection
with the centre ${\cal Z}({\frak g})$ iff $\mu$ is not a zero of an 
exotic factors of some PRV determinant.

\begin{rem}{} Taking into account that 
any completely reducible Lie superalgebra is a direct sum
of simple Lie algebras and  Lie superalgebras of the type
$\osp(1,2l)$, our results (theorems~\ref{thmhess},~\ref{thmprv},~\ref{main})
can be directly extended to the case of any 
finite dimensional completely reducible Lie 
superalgebra ${\frak g}$. 
\end{rem}

{\em Acknowledgement.} We are greatly indebted to our teacher 
Prof. A.~Joseph who proposed the problem and pushed us to solve it. 
His paper "Sur l'annulateur d'un module de Verma" was the 
main inspiration of the present work. 
We would like to thank  V.~Hinich for numerous suggestions and support.
We are also grateful to T.~Levasseur and M.~Duflo for helpful discussions
and important remarks.

\section{background}
\label{notation}
In this Section we fix the main notations we use throughout this paper.
The notation ${\Bbb N}^+$ will stand for the set of positive integers.
The base field we are going to work with is ${\Bbb C}$. 
\subsection{}
Let ${\frak g}$ be the Lie superalgebra $\osp(1,2l),\ l\geq 1$ 
(see Kac~\cite{k2} for a
presentation of this Lie superalgebra by generators and relations). Denote
by ${\frak g}_0$ the even part and by ${\frak g}_1$ 
the odd part of ${\frak g}$.
We recall that ${\frak g}_0\simeq \mbox{sp}(2l)$. Fix a Cartan subalgebra
${\frak h}$ in ${\frak g}_0$. Denote by $\Delta_0$ (resp., $\Delta_1$)
the set of even (resp., odd) roots of ${\frak g}$. 
Set $\Delta=\Delta_0\cup\Delta_1$.
Let $\Delta_{irr}$ be the set of irreducible roots 
of $\Delta$. Then $\Delta_{irr}=\overline{\Delta_0}\cup \Delta_1$, where
$\overline{\Delta_0}\!\!:=\!\Delta_0\backslash 2\Delta_1$.

Fix a basis of simple roots $\pi$ of $\Delta$, and define correspondingly 
the sets 
$\Delta^{\pm}, \Delta_0^{\pm},\Delta_1^{\pm},\overline{\Delta_0}^{\pm},
\Delta_{irr}^{\pm}$. Denote by $W$ the Weyl group of $\Delta$.
Set 
$$\begin{array}{ccc}
\rho_0\!:={1\over 2}\displaystyle\mathop{\sum}_{\alpha\in \Delta_0^+}\alpha,& 
\rho_1\!:={1\over 2}\displaystyle\mathop{\sum}_{\alpha\in \Delta_1^+}\alpha,& 
\rho\!:=\rho_0-\rho_1={1\over 2}\displaystyle\mathop{\sum}_{\alpha\in
\Delta_{irr}^+}\alpha.
\end{array}$$ 
Introduce the standard partial order relation on ${\frak h}^*$: 
$\lambda\leq \mu \Longleftrightarrow \mu-\lambda\in {\Bbb N}\pi$.

Denote by $(-,-)$ the non-degenerate bilinear form on ${\frak h}^*$ coming 
from the
restriction of the Killing form of ${\frak g}_0$
to ${\frak h}$. Let $\varphi:{\frak h}^*\longrightarrow {\frak h}$ be the 
isomorphism 
given by $\varphi(\lambda)(\mu):=(\lambda,\mu)$. For any $\lambda,\mu\in 
{\frak h}^*,\
\mu\not=0$ one defines 
$$\langle\lambda,\mu\rangle:=2\displaystyle{(\lambda,\mu)\over (\mu,\mu)}.$$

\subsection{}
\label{realization}
One has the following useful realization of $\Delta$. 
Identify ${\frak h}^*$ with ${\Bbb C}^l$ and consider $(-,-)$ as 
a scalar product on ${\Bbb C}^l$. Then there exists an orthonormal basis 
$\{\beta_1,\ldots,\beta_l\}$ such that
$$\begin{array}{lr}\pi=\{\beta_1-\beta_2,\ldots,\beta_{l-1}-\beta_l,\beta_l\}& 
\\
\Delta_0^+=\{\beta_i\pm\beta_j, 1\leq i<j\leq l,\ 2\beta_i,\ 1\leq i\leq l\},&
\Delta_1^+=\{\beta_i, 1\leq i\leq l\}\\
\Delta_{irr}^+=\{\beta_i\pm\beta_j,\ \beta_i,\ 1\leq i<j\leq l\},&
\overline{\Delta_0}^+=\{\beta_i\pm\beta_j,\ 1\leq i<j\leq l\}\\
\rho=\displaystyle\sum_{i=1}^l(l-i+{1\over 2})\beta_i,&
\rho_0=\displaystyle\mathop{\sum}_{i=1}^l(l-i+1)\beta_i
\end{array}$$
and the Weyl group $W$ is just the group of the signed permutations of the 
$\beta_i$. 
Define the translated action of $W$ on ${\frak h}^*$ by the formula: 
$w.\lambda=w(\lambda+\rho)-\rho$ for $\lambda\in {\frak h}^*$ and $w\in W$. 
For any element $w\in W$ set $\sgn w:={(-1)}^{l(w)}$, where
$l(w)$ is the length of a reduced decomposition of $w$. 

There exists a Chevalley antiautomorphism $\sigma$ of ${\frak g}$
of order $4$ which leaves invariant the elements of ${\frak h}$ and
maps ${\frak g}_{\alpha}$ to ${\frak g}_{-\alpha}$ for any $\alpha\in \Delta$.

\subsection{Enveloping algebra} 
As usual, if ${\frak k}$ is a Lie superalgebra, 
${\cal U}({\frak k})$ denotes its enveloping algebra.
Set ${\cal F}$ the natural 
filtration of ${\cal U}({\frak g})$ defined by 
${\cal F}={({\frak g}^n)}_{n\in \Bbb N}$.
The graded algebra of ${\cal U}({\frak g})$ associated to ${\cal F}$ is
the symmetric superalgebra denoted by 
${\cal S}({\frak g})\simeq {\cal S}({\frak g}_0)\otimes\bigwedge{\frak g}_1$ 
which is not a domain. Nevertheless, 
Aubry and Lemaire showed (\cite{al}) that ${\cal U}({\frak g})$ is a domain.

We define the supercentre to be the vector subspace of ${\cal U}({\frak g})$ 
generated by the homogenous elements $a$ such that $ax={(-1)}^{|a||x|}xa$
for all homogenous elements $x$ in ${\cal U}({\frak g})$. 
For ${\frak g}=\osp (1,2l)$, the supercentre coincides with the 
genuine centre.

The Lie superalgebra ${\frak g}$ acts on ${\cal U}({\frak g})$ and 
${\cal S}({\frak g})$ by superderivation
via the adjoint action. We denote these actions by $\ad$. 
Throughout this paper, 
an action of any element of ${\frak g}$ on ${\cal U}({\frak g})$ means always
the adjoint action.

We identify ${\cal U}({\frak h})$ with ${\cal S}({\frak h})$.

\subsection{Verma and simple highest weight modules}
For any $\alpha\in \Delta$, let ${\frak g}_{\alpha}$ be the one-dimensional 
subspace of weight $\alpha$ of ${\frak g}$ and
${\frak n}^{\pm}=\displaystyle\mathop{\bigoplus}_{\alpha\in \Delta^{\pm}}
{\frak g}_{\alpha}$,\ ${\frak b}^{\pm}=\frak h\oplus{\frak n}^{\pm}$.  

For a fixed ${\lambda}\in {\frak h}^*$, let $\widetilde{\Bbb C}_{\lambda}$ be 
the one dimensional ${\frak b}^+$-module with ${\frak n}^+v=0$ and
$hv=\lambda(h)v$ for all $h\in {\frak h}$ and 
$v\in \widetilde{\Bbb C}_{\lambda}$. Set
$$\widetilde{M}(\lambda)={\cal U}({\frak g}){\otimes}_{{\cal U}({\frak b}^+)}
\widetilde{\Bbb C}_{\lambda}.$$
The Verma module $\widetilde{M}(\lambda)$ has a unique simple quotient 
denoted by
$\widetilde{V}(\lambda)$.

\subsection{The $\widetilde{\cal O}$ category}
\label{ocat}
Let $M$ be a ${\frak g}$-module. For any $\lambda\in{\frak h}^*$, set 
$$M_{\lambda}=\left\{m\in M|\ hm=\lambda(h)m,\ 
\forall h\in {\frak h}\right\}.$$
A non-zero vector $v\in M$ has weight $\lambda$ if $v\in M_{\lambda}$.
For any subspace
$N$ of $M$ we denote by $\Omega(N)$ the set of weights $\lambda\in{\frak h}^*$
 such
that $N\cap M_{\lambda}\not=\{0\}$. The module $M$ is said to be 
diagonalizable if
$M=\displaystyle\mathop{\bigoplus}_{\lambda\in{\frak h}^*}M_{\lambda}$.
If $M$ is a diagonalizable module and
$\dim M_{\lambda}<\infty$ for all $\lambda\in {\frak h}^*$, we set
$\ch M=\displaystyle\mathop{\sum}_{\lambda\in{\frak h}^*}
\dim M_{\lambda}e^{\lambda}$.

Throughout this article, we shall work essentially with two categories of 
${\frak g}$-modules.
The first one is the category of completely reducible ${\frak g}$-modules.
If $M$ is such a module and 
$V$ is any simple module, we shall denote by $[M:V]$ the number 
$\dim \Hom_{\frak g} (V,M)$.
The second one, denoted by $\widetilde{\cal O}$, is the full 
subcategory of the category of 
${\frak g}$-modules whose objects are ${\frak g}$-modules $M$ such that 
\begin{enumerate}
\item $M$ is diagonalizable
\item $\forall \lambda\in \Omega(M),\dim M_{\lambda}<+\infty$
\item $\exists \lambda_1,\ldots,\lambda_k\in {\frak h}^* |\ 
\Omega(M)\subset \mathop{\bigcup}_{i=1,\ldots, k}(\lambda_i-{\Bbb N}\pi)$
\end{enumerate}

Any object of the $\widetilde{\cal O}$ category, considered as
a ${\frak g}_0$-module, belongs to the ${\cal O}$ category 
(see~\cite{jfr}, 3.5 for definition). 
In particular, any ${\frak g}$-module $M$ of $\widetilde{\cal O}$
has finite length. If $V$ is a simple highest weight module, we denote 
by $[M:V]$ the number of simple quotients isomorphic to $V$ in any 
composition series of $M$. This number does not depend of the choice of the 
composition series of $M$.

\subsection{Finite dimensional representations}
\label{findim}
Define for $r\in\{1,\ldots,l\}$ the fundamental weight
$\omega_r=\displaystyle\sum_{i=1}^r\beta_i$,  and introduce the set  
$$P^+(\pi)\!:=\displaystyle\sum_{r=1}^l{\Bbb N}\omega_r=
\{\lambda\in {\frak h}^*|\ \langle\lambda,\beta_l\rangle\in 2{\Bbb N},\
\langle\lambda,\beta_i-\beta_{i+1}\rangle\in {\Bbb N},\ 
\forall i=1,\ldots,l-1\}.$$
Kac (see~\cite{k1}) showed that  $\widetilde{V}(\lambda)$ is 
finite dimensional iff
$\lambda\in P^+(\pi)$. For any $\mu\in{\frak h}^*$, let $V(\mu)$ be the simple
${\frak g}_0$-module of highest weight $\mu$. Remark that
$\{\beta_1-\beta_2,\ldots,\beta_{l-1}-\beta_l,2\beta_l\}$ is a basis of 
simple roots of $\Delta_0$  and that $\langle \mu,2\beta_l\rangle={1\over
2}\langle\mu,\beta_l\rangle$ for all $\mu\in {\frak h}^*$.
Thus $V(\lambda)$ is finite dimensional iff $\lambda\in P^+(\pi)$.

\subsubsection{}
\label{duality} Observe that the longest element of the Weyl group acts
on ${\frak h}^*$ as $-id$. Therefore for any finite dimensional module $V$,
$\Omega(V)=-\Omega(V)$. 

\section{The separation theorem}
\label{linkann}
Recall the separation theorem established by Musson in~\cite{mu},1.4:
\subsection{}
\begin{thm}{separation} There exists an ad-invariant subspace ${\cal H}$ in 
${\cal U}({\frak g})$ such that the multiplication map induces an
$\ad {\frak g}$-isomorphism ${\cal U}({\frak g})\simeq
{\cal Z}({\frak g})\otimes {\cal H}$. Moreover, for every simple
finite dimensional module $\widetilde V$, $[{\cal H}:\widetilde{V}]=\dim 
\widetilde{V}_0$.
\end{thm}

\subsection{}
\label{redtosimple}
Since the centre ${\cal Z}({\frak g})$ acts by a scalar on  every Verma module 
$\widetilde{M}(\mu)$, the separation theorem implies
the following equivalence :
$$Ann_{{\cal U}({\frak g})}\widetilde{M}(\mu)={\cal U}({\frak g})
\Ann_{{\cal Z}({\frak g})}\widetilde{M}(\mu)\ \Longleftrightarrow \ 
\Ann_{\cal H}\widetilde{M}(\mu)=\{0\}.$$
As in the classical case, any Verma module contains a simple Verma submodule.
Let $\widetilde{M}(\mu)$ be a Verma module and $M$ be its simple
Verma submodule. Obviously, $$\Ann_{{\cal H}}\widetilde{M}(\mu)\subset 
\Ann_{{\cal H}} M.$$
Since $\widetilde{M}(\mu)$ and $M$ have the same central character,
$M$ has the form $M\simeq M(w.\mu)$ for some $w\in W$ (see~\cite{mu},1.1). 
By definition of the 
translated action of $W$,  $w.\mu+\rho=w(\mu+\rho)$ and as $W$ acts on 
$\Delta_1$ by signed permutations, one has
$$\exists \alpha\in \Delta_1^+,\ (\mu+\rho,\alpha)=0 \Longleftrightarrow 
\exists \alpha\in \Delta_1^+,\ (w.\mu+\rho,\alpha)=0.$$

Hence the proof of the implication "$\Leftarrow$" of our
main~\Thm{main} 
can be reduced to the case of simple Verma modules. Indeed for this, 
it is enough to
show that
\begin{equation}
\label{eqiv1} 
\bigr[\widetilde{M}(\mu)\mbox{ simple and } (\mu+\rho,\alpha)\not=0\ \forall 
\alpha\in \Delta_1^+\bigr] \ \Longrightarrow \ \ \Ann_{{\cal H}}\widetilde{M}
(\mu)=0.
\end{equation}

\section{Factorization of the Shapovalov determinants}
\label{shapo}
\subsection{}
\label{contrform}
The triangular decomposition ${\cal U}({\frak g})={\cal U}({\frak n}^-)\otimes 
{\cal U}({\frak h})\otimes
{\cal U}({\frak n}^+)$ leads to the following decomposition 
$${\cal U}({\frak g})={\cal U}({\frak h})\oplus({\frak n}^-{\cal U}({\frak g})+
{\cal U}({\frak g}){\frak n}^+)$$
The projection $\Upsilon:{\cal U}({\frak g})\longrightarrow 
{\cal U}({\frak h})$ 
with respect to this decomposition is called the
Harish-Chandra projection. Observe that $\Upsilon(ha)=h\Upsilon(a)$ for
all $h\in {\cal U}({\frak h}),a\in {\cal U}({\frak g})$.

The Harish-Chandra projection map allows us to define a contravariant form on 
each Verma 
module $\widetilde{M}(\lambda)$ by 
the formula $${\langle av_{\lambda},bv_{\lambda}\rangle}_{\lambda}=
\Upsilon(\sigma(a)b)
(\lambda)\quad \forall a,b\in {\cal U}({\frak g}),$$
where $v_{\lambda}$ is a primitive vector of highest weight $\lambda$. 

The kernel of this form is the largest proper submodule of $\widetilde{M}
(\lambda)$.
Consequently, the Verma 
module $\widetilde{M}(\lambda)$ is simple if and only if the form 
${\langle-,-\rangle}_{\lambda}$ is non-degenerate.
The subspaces of different weights are pairwise orthogonal with respect to 
${\langle-,-\rangle}_{\lambda}$. Therefore, the quest 
for a criterion of simplicity for a Verma module leads naturally to the
 analysis of the zeroes of the 
so called Shapovalov determinants: $\det S_{\nu}$, where $\nu\in {\Bbb N}\pi$ 
and
$S_{\nu}:\ {{\cal U}({\frak n}_-)}_{\nu}\times {{\cal U}({\frak n}_-)}_{\nu}
\mapsto {\cal S}({\frak h})$ is
defined by the formula $$S_{\nu}(x,y)=\Upsilon(\sigma(x)y).$$

The factorization of these determinants was established by Kac 
(see~\cite{k3}) for all classical simple Lie superalgebras.
For the present case ${\frak g}=\osp(1,2l)$ (see also~\cite{mu2} 2.4)
 one has the

\subsection{}
\begin{thm}{shapovalov} For all $\nu\in {\Bbb N}$, one has, up to a non-zero
scalar,

$\vtop{\hbox to \textwidth{$\det {S}_{\nu}=\displaystyle
\mathop{\prod}_{\alpha\in\overline{\Delta_0^+}}\mathop{\prod}_{m=1}^{\infty}
{(\varphi(\alpha)+(\rho,\alpha)-{1\over 2}m(\alpha,\alpha))}^
{\tau(\nu-m\alpha)}
\times$\hfill}
\hbox to \textwidth{\hfill$\displaystyle\mathop{\prod}_{\alpha\in{\Delta_1^+}}
\mathop{\prod}_{m=1}^{\infty}
{(\varphi(\alpha)+(\rho,\alpha)-{1\over 2}(2m-1)(\alpha,\alpha))}^
{\tau(\nu-(2m-1)\alpha)}$}}$
where $\tau:{\Bbb Z}\pi\longrightarrow {\Bbb N}$ is the Kostant 
partition function defined by $\tau({\Bbb Z}\pi\backslash {\Bbb N}\pi)=0$
and 
$$\tau(\nu)=\# \left\{ {\{k_{\alpha}\}}_{\alpha\in \Delta^+} |\ 
  k_{\alpha}\in {\Bbb N}\mbox{ and }
\sum_{\alpha\in\Delta^+}k_{\alpha}\alpha=\nu\right\}, \ \ \forall
 \nu\in{\Bbb N}\pi.$$
\end{thm}

\subsection{} 
\label{weightsub}
Take $\lambda\in {\frak h}^*$.
The theorem above allows us to describe the weights of the largest proper
submodule $N$ of the Verma module $\widetilde{M}(\lambda)$:
$$\begin{array}{c}
\Omega (N)=\displaystyle\mathop{\bigcup}_{\alpha\in \Delta_{\lambda}}
\{\lambda-\langle\lambda+\rho,\alpha\rangle\alpha-{\Bbb N}\pi\}, 
\text { where } \\
\Delta_{\lambda}\:\!=\{\alpha\in \overline{\Delta^+_0}|\ 
\langle\lambda+\rho,\alpha\rangle\in {\Bbb N}^+\}\cup \{\beta\in \Delta^+_1|\ 
\langle\lambda+\rho,\beta\rangle\in 2{\Bbb N}+1\}.
\end{array}$$
In particular, it implies the following criterion of simplicity for 
$\widetilde{M}(\lambda)$:
\subsection{}
\begin{cor}{simpleV}
The Verma module $\widetilde{M}(\lambda):\ \lambda\in {\frak h}^*$ is 
simple if and only if :
$$\begin{array}{lr} \langle\lambda+\rho,\alpha\rangle\not\in{\Bbb N}^+,&
\forall \alpha\in \overline{\Delta_0^+}\\
 \langle\lambda+\rho,\alpha\rangle\not\in(2{\Bbb N}+1),& \forall \alpha\in
{\Delta_1^+}.
\end{array}$$
\end{cor}

\section{Hesselink Formula}
\label{hesselink}
The goal of this section is to give a formula describing a multiplicity
of every simple finite dimensional module $\widetilde {V}(\lambda)$ 
in the graded component $H^n$ where $H:=\gr_{\cal F}{\cal H}$.
More precisely, let ${\cal S}^n({\frak g})
\subseteq {\cal S}({\frak g})$ be the subspace
of homogeneous elements of degree $n$. For a finite dimensional module 
$\widetilde {V}(\lambda)$ we will give an explicit formula for 
a Poincar\'e series $P_{\lambda}(q)$ of $H=\gr_{\cal F} {\cal H}$
which is defined by
$$P_{\lambda}(q)\!:=\sum_{n=0}^{\infty} [H^n:{\widetilde V}(\lambda)]q^n,
\text { where } H^n\!:={\cal S}^n({\frak g })\cap H.$$

\subsection{Notation}
\label{nthess}
Consider the group ring ${\Bbb Z}[{\frak h}^*]$. Its elements
are finite sums $\sum_{\lambda\in {\frak h}^*} c_{\lambda}e^{\lambda}, \
c_{\lambda}\in {\Bbb Z}$.
Our characters lie in the  power series ring ${\Bbb Z}[{\frak h}^*][[q]]$.
For each $\mu\in{\frak h}^*$ define the ${\Bbb Z}[[q]]$-linear 
homomorphism $\pi_{\mu}$ by 
$$\pi_{\mu}: {\Bbb Z}[{\frak h}^*][[q]]\ \rightarrow\ {\Bbb Z}[[q]],\ \ \ \ 
\ \ e^{\lambda}\mapsto \delta_{\lambda, \mu},
\ \forall \lambda\in {\frak h}^*.$$

\subsubsection{}
\label{jprop}
For each $w\in W$ define the automorphism of the ring  
${\Bbb Z}[{\frak h}^*][[q]]$ by 
$$w(e^{\lambda})\!:=e^{w\lambda},\ w(q)\!:=q.$$
Let $J$ be the ${\Bbb Z}[[q]]$-linear endomorphism of 
${\Bbb Z}[{\frak h}^*][[q]]$ given by the formula
$$J\!:=\sum_{w\in W} \sgn (w)w.$$

We use the following properties of the operator $J$:
\begin{enumerate}
\item $J(wa)=\sgn (w)J(a)$ for all $w\in W$. Consequently, if there exists $w$
such that $wa=a$ and  $\sgn (w)=-1$ then $J(a)=0$.
\item Since for any $\lambda\in P^+(\pi)$ the stabilizer
$\Stab_W\lambda$ is generated by the simple reflections it 
contains~(\cite{jq}, A.1.1) and $P(\pi)=WP^+(\pi)$, one
can conclude from (i) that 
$$\Stab_W\mu\not=\{\id\}\ \Longrightarrow J(e^{\mu})=0$$
for any $\mu\in P(\pi)$.
\item Fix a subgroup $K$ of $W$ and representatives 
$g_1,\ldots, g_r$ of the left cosets $W/K$. Set
$$J^{K}\!:=\sum_{w\in K} \sgn(w) w,\ \ \ \ 
J^{W/K}\!:=\sum_{i=1}^r \sgn(g_i){g_i}.$$
Then $J^{W/K}J^{K}=J$.
\item Assume that $wa=a$ for all $w\in K$. Then $J^K(ab)=aJ^K(b)$.
\end{enumerate}

\subsubsection{}
\label{sprth}
For a graded ${\frak h}$-submodule $N$ of ${\cal S}({\frak g})$ 
its graded character $\ch_q N$ is defined to be the element 
of ${\Bbb Z}[{\frak h}^*][[q]]$ given by the formula
$$\ch_q N\!:=\sum_{r=0}^{\infty} \ch (N\cap {\cal S}^r({\frak g})) q^r.$$
The following relations hold
$$\begin{array}{l}
\ch_q\Lambda {\frak g}_1=\displaystyle\sum_{r=0}^{2l} 
(\ch\Lambda^r {\frak g}_1)q^r=\displaystyle\prod_{\beta\in\Delta_1}
(1+qe^{\beta})\\
\ch_q{\cal S}({\frak n}_0^{\pm})=\displaystyle\prod_{\alpha\in\Delta_0^+}
(1-qe^{\pm\alpha})^{-1},\\
\ch_q{\cal S}({\frak n}^{\pm})=\displaystyle\prod_{\alpha\in\Delta^+_0}
(1-qe^{\pm\alpha})^{-1}\!\displaystyle\prod_{\beta\in\Delta^+_1}
(1+qe^{\pm\beta}),\\
\ch_q{\cal S}({\frak n}^+\oplus {\frak n}^-)=
\ch_q{\cal S}({\frak n}_0^+\oplus {\frak n}_0^-)\cdot
\ch_q\Lambda {\frak g}_1.
\end{array}
$$

Moreover, the separation theorem implies that 
$$\ch_q {\cal S}({\frak g})=\ch_q {{\cal S}({\frak g})}^{{\frak g}}\ch_q H.$$
By~\cite{mu}, 1.2 ${{\cal S}({\frak g})}^{{\frak g}}\cong 
{{\cal S}({\frak h})}^W$.
By~\cite{jfr}, 8.7
$$\ch_q {\cal S}({\frak h})=\ch_q {\cal S}({\frak h})^W\!\displaystyle
\sum_{w\in W}q^{l(w)}$$
so
$$\ch_q {\cal S}({\frak g})=\ch_q {\cal S}({\frak n}^-\oplus{\frak n}^+)
\ch_q {\cal S}({\frak h})=
\ch_q {\cal S}({\frak n}^-\oplus{\frak n}^+)
\ch_q {{\cal S}({\frak g})}^{{\frak g}}
\!\displaystyle\sum_{w\in W}q^{l(w)}.$$
Hence
\begin{equation}\label{hch}
\ch_q H=\ch_q {\cal S}({\frak n}^-\oplus{\frak n}^+)\!\displaystyle
\sum_{w\in W}q^{l(w)}=\prod_{\alpha\in\Delta_0}(1-qe^{\alpha})^{-1}
\ch_q\Lambda {\frak g}_1\!\!\sum_{w\in W}q^{l(w)}.
\end{equation}

\subsection{Classical Case}
\label{hsscl}
We use the previous notations for the corresponding classical objects.
For the classical case the Hesselink formula (see~\cite{jfr}, 9.10) states that
$P_{\lambda}(q)$ is equal to the coefficient of $e^0=1$ in the 
expression
\begin{equation}
\label{pseri}
Q_{\lambda}(q)\!:=e^{-\rho}J(e^{\lambda+\rho})\prod_{\alpha\in\Delta^+}
(1-qe^{-\alpha})^{-1}.
\end{equation}
The proof (see~\cite{jfr}, 8.6) is essentially based on the equality 
\begin{equation}
\label{nicefrm}
J(e^{\rho}\!\prod_{\alpha\in\Delta^+}(1-qe^{-\alpha}))=J(e^{\rho})\!\!
\sum_{w\in W} q^{l(w)}.
\end{equation}

\subsection{}
A key point of the proof of~\Thm{thmhess} is the following 
analogue of equality~\ref{nicefrm}.
\begin{prop}{hesprop}
\begin{equation}
\label{hesfr}
J\left(e^{\rho}\!\biggl(\prod_{\alpha\in{\Delta^+_0}}(1-qe^{-\alpha})\biggr)\!
(1+q^{2l}e^{\beta_1})
\!\!\prod_{\beta\in {\Delta^+_1\setminus\{\beta_1\}}}\!\!(1+e^{\beta})\right)
=J(e^{\rho})\ch_q\Lambda {\frak g}_1\!\!\sum_{w\in W}q^{l(w)}
\end{equation}
where $\beta_1$ is the maximal odd root.
\end{prop}
\begin{pf}
The proposition will be proven in~\ref{heslem}---\ref{endpfprop}.
First of all we need the following technical lemma.
\subsubsection{}
\begin{lem}{heslem}
Set $\Gamma_{\bullet}:=\{\displaystyle\sum_{i=2}^l k_i\beta_i,\ \ 
k_i\in\{0,1\}\}$.
For any $r=0,1,\ldots, 2l-1$ there exists a unique
$\gamma_r\in \Gamma_{\bullet}$ such that 
$\Stab_W(\rho-r\beta_1+\gamma_r)=\{\id\}$. Moreover 
$\rho-r\beta_1+\gamma_r=w_r\rho$ for some $w_r\in W$ and $\sgn (w_r)=(-1)^r$.
\end{lem}
\begin{pf}
Recall that the set $\{\beta_i\}_{i=1}^l$ is an orthonormal basis of
${\frak h}^*$ and $W$ acts on this basis by signed permutations.
Therefore the stabilizer of $\lambda\in {\Bbb Q}\pi$ is trivial
iff $|(\lambda,\beta_i)|$ are pairwise distinct non-zero values.

Assume that $\Stab_W(\rho-r\beta_1+\gamma)=\{\id\}$
for some $\gamma\in \Gamma_{\bullet}$ and set $\lambda:=\rho-r\beta_1+\gamma$.
Recall that $(\rho,\beta_i)=l+1/2-i$. Since $0\leq r<2l$ one has
$$\{|(\lambda,\beta_i)|\}_{i=1}^l\subseteq \{l+1/2-i\}_{i=1}^l
=\{|(\rho,\beta_i)|\}_{i=1}^l.$$
Taking into account that all $|(\lambda,\beta_i)|$ are pairwise distinct 
we conclude that 
$\{|(\lambda,\beta_i)|\}_{i=1}^l=\{|(\rho,\beta_i)|\}_{i=1}^l$.
This implies that $\lambda=w\rho$ for some $w\in W$ and that $\gamma$
is determined by the value of $(\lambda,\beta_1)=l-1/2-r$. 
The explicit expressions  of $\gamma$ and $w$ for a fixed $r$ are:
$$
\left\{
\begin{array}{lll}
r=0\  & \gamma=0, & w=id \\
1\leq r<l\  & \gamma=\displaystyle\sum_{i=2}^{r+1}\beta_i, & w=s_1\ldots s_r \\
l\leq r<2l\  & \gamma=\displaystyle\sum_{i=2}^{2l-r}\beta_i, & w=s_1\ldots s_l
s_{l-1}\ldots s_{2l-r} 
\end{array}
\right.
$$
This completes the proof of the lemma.
\end{pf}
\subsubsection{}
Let $\Delta_{\bullet}$ be the subsystem of $\Delta_0$ generated by
$\pi\setminus\{\beta_1-\beta_2\}$ and $W_{\bullet}$ be the corresponding Weyl 
group that is the subgroup of $W$ generated by the reflections
$s_{\alpha_2},\ldots,s_{\alpha_l}$.
Set $\Delta_{\bullet}^+:=\Delta_{\bullet}\cap \Delta^+$,  
$\rho_{\bullet}:=\displaystyle\sum_{\alpha\in \Delta_{\bullet}^+}\alpha/2.$
Observe that $\Delta_{\bullet}$ is the root system corresponding
to the simple Lie algebra $\ssp (l-1)$ and that 
$\rho_{\bullet}=\rho_0-l\beta_1$. By~\cite{h}, 3.15 for
the Weyl group $W_n$ of the Lie algebra $\ssp(n)$ one has
$\displaystyle\sum_{w\in W_n}q^{l(w)}=
(1-q)^{-n}\displaystyle\prod_{i=1}^n (1-q^{2i})$
and thus
\begin{equation}
\label{qW}
(1-q^{2l})\!\!\sum_{w\in W_{\bullet}} q^{l(w)}=(1-q)\!\!\sum_{w\in W}q^{l(w)}
\end{equation}
\subsubsection{}
\label{endpfprop}
Observe that
\begin{equation}
\label{ex}
\begin{array}{lcl}
(1-q)\displaystyle\prod_{\alpha\in {\Delta_{0}^+}\setminus{\Delta_{\bullet}^+}}
(1-qe^{-\alpha})&=&
(1-qe^{-2\beta_1})(1-q)\displaystyle\prod_{i=2}^l 
(1-qe^{-\beta_1-\beta_i})(1-qe^{-\beta_1+\beta_i})\\
&=&\displaystyle\prod_{i=1}^l 
(1-qe^{-\beta_1-\beta_i})(1-qe^{-\beta_1+\beta_i})\\
&=&\displaystyle\sum_{r=0}^{2l}\ch (\Lambda^r {\frak g}_1)\!\cdot\!
(-q)^re^{-r\beta_1}.
\end{array}
\end{equation}
The characters $\ch (\Lambda^r {\frak g}_1)$ are stable under
the action of $W$ since ${\frak g}_1$ is a ${\frak g}_0$-module.
By the property (iv) of~\ref{jprop} with $K=W$, one has
\begin{equation}
\label{cc0}
\begin{array}{lcr}
\lefteqn{
(1-q)J\left(e^{\rho}\!\biggl(\displaystyle\prod_{\alpha\in{\Delta^+_0}}
(1-qe^{-\alpha})\biggr)\!
(1+q^{2l}e^{\beta_1})
\!\!\displaystyle\prod_{\beta\in {\Delta^+_1\setminus\{\beta_1\}}}
\!(1+e^{\beta})\right)=}\\
& &\displaystyle\sum_{r=0}^{2l}\ch (\Lambda^r {\frak g}_1)
(-q)^rJ\biggl(e^{\rho-r\beta_1}(1+q^{2l}e^{\beta_1})
\prod_{\alpha\in {\Delta_{\bullet}^+}}(1-qe^{-\alpha})
\prod_{i=2}^l (1+e^{\beta_i})\biggr).
\end{array}
\end{equation}

Decompose $J=J^{W/W_{\bullet}}J^{W_{\bullet}}$ (see (iii) of~\ref{jprop}) and
observe that the elements 
$$e^{\beta_1} \text { and }
e^{\rho-\rho_{\bullet}}\displaystyle\prod_{i=2}^l (1+e^{\beta_i})=
e^{l\beta_1}\displaystyle\prod_{i=2}^l (e^{-\beta_i/2}+e^{\beta_i/2})$$
are stable under the action of $W_{\bullet}$.
Using the formula~(\ref{nicefrm}) with respect to $W_{\bullet}$ one can 
rewrite the right hand side of~(\ref{cc0}) as
\begin{equation}
\label{cc}
\begin{array}{l}
\displaystyle\sum_{r=0}^{2l}\ch (\Lambda^r {\frak g}_1)
(-q)^rJ^{W/W_{\bullet}}\biggl(e^{\rho-\rho_{\bullet}-r\beta_1}(1+q^{2l}
e^{\beta_1})
\prod_{i=2}^l (1+e^{\beta_i}) J^{W_{\bullet}}\Bigl(e^{\rho_{\bullet}}
\prod_{\alpha\in{\Delta_{\bullet}^+}}(1-qe^{-\alpha})\Bigr)\biggr)=\\
\displaystyle\sum_{r=0}^{2l}\ch (\Lambda^r {\frak g}_1)
(-q)^rJ^{W/W_{\bullet}} \biggl(e^{\rho-\rho_{\bullet}-r\beta_1}(1+q^{2l}
e^{\beta_1})
\prod_{i=2}^l (1+e^{\beta_i})J^{W_{\bullet}}(e^{\rho_{\bullet}})\!\!
\displaystyle\sum_{w\in W_{\bullet}} q^{l(w)}\biggr)=\\
\Bigl(\displaystyle\sum_{w\in W_{\bullet}} q^{l(w)}\Bigr)
\displaystyle\sum_{r=0}^{2l}\ch (\Lambda^r {\frak g}_1)(-q)^r
J\biggl(e^{\rho-r\beta_1}(1+q^{2l}e^{\beta_1})\prod_{i=2}^l 
(1+e^{\beta_i})\biggr)=\\
\Bigl(\displaystyle\sum_{w\in W_{\bullet}} q^{l(w)}\Bigr)
\displaystyle\sum_{r=0}^{2l}\ch (\Lambda^r {\frak g}_1)a_r,
\end{array}
\end{equation}
where $a_r\:=(-q)^r
J\biggl(e^{\rho-r\beta_1}(1+q^{2l}e^{\beta_1})\prod_{i=2}^l (1+e^{\beta_i})
\biggr)$.

One can simplify the expression obtained in the following way.
For $r\in\{1,\ldots,2l-1\}$ take $w_r,w_{r-1}$ as in~\Lem{heslem}. Then
\begin{equation}
\label{coeffr}
a_r=(-q)^rJ(e^{w_r\rho}+q^{2l}e^{w_{r-1}\rho})=
(1-q^{2l})q^rJ(e^{\rho}).
\end{equation}

For $r=0,2l$ we have $\ch (\Lambda^0 {\frak g}_1)=
\ch (\Lambda^{2l} {\frak g}_1)=1$,
therefore
$$a_0+a_{2l}=J\left((e^{\rho}+q^{2l}e^{\rho-2l\beta_1})(1+q^{2l}e^{\beta_1})
\displaystyle\prod_{i=2}^l (1+e^{\beta_i})\right).$$
Consider the reflection $s=s_{\beta_1}\in W$.
One has $s\rho=\rho-(2l-1)\beta_1$. Thus the expression 
$(e^{\rho+\beta_1}+e^{\rho-2l\beta_1})\displaystyle\prod_{i=2}^l 
(1+e^{\beta_i})$ is stable under the action of $s$.
Since $\sgn (s)=-1$, property (i) of~\ref{jprop} implies that
\begin{equation}
\label{cc2}
J\left((e^{\rho+\beta_1}+e^{\rho-2l\beta_1})\prod_{i=2}^l 
(1+e^{\beta_i})\right)=0.
\end{equation}

Using~\Lem{heslem} with respect to $r=0$ and $r=2l-1$, we obtain
\begin{equation}
\label{cc1}
J\left((e^{\rho}+q^{4l}e^{\rho-(2l-1)\beta_1})\displaystyle
\prod_{i=2}^l (1+e^{\beta_i})\right)=
J(e^{\rho}-q^{4l}e^{\rho})=J(e^{\rho})(1-q^{2l})(1+q^{2l}).
\end{equation}
Summarizing~(\ref{cc2}) and~(\ref{cc1}) we get
\begin{equation}
\label{coeff0}
a_0+a_{2l}=J\left((e^{\rho}+q^{2l}e^{\rho-2l\beta_1})(1+q^{2l}e^{\beta_1})
\displaystyle\prod_{i=2}^l (1+e^{\beta_i})\right)=J(e^{\rho})(1-q^{2l})
(1+q^{2l}).
\end{equation}

By the substitution of~(\ref{coeffr}) and~(\ref{coeff0}) into~(\ref{cc}) 
one can rewrite~(\ref{cc0}) as
\begin{equation}
\label{c}
\begin{array}{l}
(1-q)J\left(e^{\rho}\!\biggl(\displaystyle\prod_{\alpha\in{\Delta^+_0}}
(1-qe^{-\alpha})\biggr)\!(1+q^{2l}e^{\beta_1})
\!\!\displaystyle\prod_{\beta\in {\Delta^+_1\setminus\{\beta_1\}}}\!\!
(1+e^{\beta})\right)=\\
\displaystyle\sum_{r=0}^{2l}\ch (\Lambda^r {\frak g}_1)q^rJ(e^{\rho})(1-q^{2l})
\left(\displaystyle\sum_{w\in W_{\bullet}} q^{l(w)}\right)=\\
\displaystyle\sum_{r=0}^{2l}\ch (\Lambda^r {\frak g}_1)q^rJ(e^{\rho})(1-q)
\left(\displaystyle\sum_{w\in W} q^{l(w)}\right)\ \text { by~(\ref{qW})}
\end{array}
\end{equation}
This completes the proof of~\Prop{hesprop}.
\end{pf}
The following theorem is the analogue of Hesselink formula for $\osp(1,2l)$.
\subsection{}
\label{hss2}
\begin{thm}{thmhess}
For any $\lambda\in P^+(\pi)$, the Poincar\'e series
$P_{\lambda}(q)$ is equal to the coefficient of $e^0=1$ in the  
expression
\begin{equation}
\label{eqthmhess}
Q_{\lambda}(q)\!:=e^{-\rho}J(e^{\lambda+\rho})\!
\biggl(\prod_{\alpha\in \Delta^+_0}(1-q e^{-\alpha})^{-1}\biggr)\!
(1+q^{2l}e^{-\beta_1})\!\!\prod_{\beta\in 
{\Delta^+_1\setminus\{\beta_1\}}}\!\!
(1+e^{-\beta})
\end{equation}
where $\beta_1$ is the maximal odd root.
\end{thm}
\begin{pf}
Observe that $\{\ch {\widetilde V(\lambda)}\!:\lambda\in P^+(\pi)\}$ are 
linearly 
independent. Hence in order to find the Poincar\'e series $P_{\lambda}(q)$
one can decompose the graded character $\ch_q H$ in terms of 
$\{\ch {\widetilde V(\lambda)}\!:\lambda\in P^+(\pi)\}$.

Recall the character formula of 
$\widetilde{V}(\lambda),\ \lambda\in P^+(\pi),$ computed by 
Kac in~\cite{k1}, (3):
\begin{equation}
\label{kacch}
\ch {\widetilde V}(\lambda)=e^{-\rho}J(e^{\lambda+\rho})
\!\prod_{\alpha\in {\Delta_0^+}}(1-e^{-\alpha})^{-1}
\!\prod_{\beta\in {\Delta_1^+}}(1+e^{-\beta}).
\end{equation}
Therefore
$$J(e^{\rho})\ch_q H=J(e^{\rho})\!\displaystyle\sum_{\lambda\in P^+(\pi)}\!
P_{\lambda}(q)\ch {\widetilde V}(\lambda)=\!\sum_{\lambda\in P^+(\pi)
}P_{\lambda}(q)J(e^{\lambda+\rho}), \text { so }$$
\begin{equation}\label{piq}
P_{\lambda}(q)=\pi_{\lambda+\rho}\bigl(J(e^{\rho})\ch_q H\bigr).
\end{equation}
On the other hand,
\begin{equation}
\begin{array}{lcl}
\label{chqH}
J(e^{\rho})\ch_q H&=&J(e^{\rho})\prod_{\alpha\in\Delta_0}(1-qe^{\alpha})^{-1}
\ch_q\Lambda {\frak g}_1\sum_{w\in W}q^{l(w)} \text { by~(\ref{hch}) }\\
&=&\displaystyle\prod_{\alpha\in\Delta_0}(1-qe^{\alpha})^{-1}
J\left(e^{\rho}\!\biggl(\displaystyle\prod_{\alpha\in
{\Delta^+_0}}(1-qe^{-\alpha})\biggr)
(1+q^{2l}e^{\beta_1})\displaystyle\prod_{i=2}^l (1+e^{\beta_i})\right)\\
&& \hskip 4truecm\text{by ~\Prop{hesprop}}\\
&=&J\left(e^{\rho}\!\biggl(\displaystyle
\prod_{\alpha\in {\Delta^+_0}}(1-qe^{\alpha})^{-1}
\biggr)(1+q^{2l}e^{\beta_1})\displaystyle\prod_{i=2}^l (1+e^{\beta_i})\right)
 \\ 
&& \hskip 4truecm\text{by property (iv) of ~\ref{jprop} with $K=W$.}
\end{array}
\end{equation}

Define the ${\Bbb Z}[[q]]$-automorphism $\iota$
of the ring ${\Bbb Z}[{\frak h}^*][[q]]$ by
$$\iota: e^{\lambda}\mapsto  e^{-\lambda},\ \forall \lambda\in{\frak h}^*.$$
Clearly $\pi_0\circ\iota=\pi_0$.
For any $a\in {\Bbb Z}[{\frak h}^*][[q]]$ one has
\begin{equation}
\label{iota}
\begin{array}{lcl}
\pi_{\mu}J(a)&=&\!\displaystyle\sum_{w\in W} \sgn (w)\pi_{w\mu}(a)=\!
\displaystyle\sum_{w\in W} \sgn (w)\pi_0 (e^{-w\mu}a)=
\displaystyle\sum_{w\in W} \sgn (w)\pi_0\circ\iota (e^{-w\mu}a)\\
&=&\pi_0\left(\displaystyle\sum_{w\in W} \sgn (w)(e^{w\mu}\iota(a))\right)=
\pi_0\left(J(e^{\mu})\iota(a)\right).
\end{array}
\end{equation}

It was already shown (see~(\ref{chqH})) that
$$P_{\lambda}(q)=\pi_{\lambda+\rho} J\biggl(e^{\rho}
\displaystyle\prod_{\alpha\in {\Delta^+_0}}(1-qe^{\alpha})^{-1}
(1+q^{2l}e^{\beta_1})\displaystyle\prod_{i=2}^l (1+e^{\beta_i})\biggr).$$
Using~(\ref{iota}) we get
$$P_{\lambda}(q)=\pi_0 \biggl(J(e^{\lambda+\rho})e^{-\rho}
\displaystyle\prod_{\alpha\in {\Delta^+_0}}(1-qe^{-\alpha})^{-1}
(1+q^{2l}e^{-\beta_1})\displaystyle\prod_{i=2}^l (1+e^{-\beta_i})\biggr).$$
This establishes the theorem.\end{pf}

\subsubsection{}
\begin{rem}{remq1}
Observe that $Q_{\lambda}(1)=\ch {\widetilde V}(\lambda)$
(see~(\ref{kacch})). Therefore $P_{\lambda}(1)$ is
equal to the coefficient of $e^0=1$ in $\ch {\widetilde V}(\lambda)$ that
is $\dim {\widetilde V}(\lambda)_0$.
On the other hand, one has
$$P_{\lambda}(1)=\sum_{n=0}^{\infty} [H^n:{\widetilde V}(\lambda)]=
[H:{\widetilde V}(\lambda)]=[{\cal H}:{\widetilde V}(\lambda)]$$
that again gives the second assertion
of the separation theorem 
$[{\cal H}:{\widetilde V}(\lambda)]=\dim {\widetilde V}(\lambda)_0$.
\end{rem}

\subsubsection{}
\begin{rem}{hesrem}
One has
$$\pi_{\lambda+\rho}J(a)=\!\displaystyle\sum_{w\in W} \sgn (w)
\pi_{w(\lambda+\rho)}(a)=\!\displaystyle\sum_{w\in W} \sgn (w)
\pi_{w.\lambda}(e^{-\rho}a)$$
where $w.\lambda=w(\lambda+\rho)-\rho$ is a twisted action
of $W$ on ${\frak h}^*$. Thus one can reformulate~\Thm{thmhess} 
as follows: for all $\lambda\in P^+(\pi),\ n\in {\Bbb N}$

$$[H^n:\widetilde {V}(\lambda)]=\sum_{w\in W} (-1)^{l(w)} P_n(w.\lambda),$$
where coefficients $P_n(\mu)$ are zero for 
$\mu\in  P(\pi)\!\setminus\!{\Bbb N}\pi$
and for $\mu\in {\Bbb N}\pi$ are given by the formula
$$
\biggl(\prod_{\alpha\in \Delta^+_0}(1-q e^{\alpha})^{-1}\biggr)
(1+q^{2l}e^{\beta_1})\prod_{\beta\in {\Delta^+_1\setminus\{\beta_1\}}}
(1+e^{\beta})=
\sum_{r=0}^{\infty}\sum_{\nu\in {\Bbb N}\pi}P_r(\nu) e^{\nu} q^r.
$$
\end{rem}

\section {The Parthasarathy--Ranga-Rao--Varadarajan determinants}
\label{prv}
The main result of this Section--~\Thm{thmprv}-- calculates the factorization 
of the Parthasarathy--Ranga-Rao--Varadarajan (PRV) determinants. 
In~\ref{prvlem},~\ref{prvcons} 
we give a construction of these determinants, in~\ref{l1} we study in 
details the case $l=1$, and
in~\ref{estimexotic} we compute a bound for the total degrees of the 
PRV determinants. The factorization formula for the PRV determinants 
is given in~\ref{thmprv}. Subsections~\ref{stfactor}--~\ref{endpfprv} 
are devoted to the proof of this formula.
\subsection{}
\label{prvlem}
In order to introduce the Parthasarathy--Ranga-Rao--Varadarajan (PRV) 
determinants, start with a little lemma. 
Retain notation of~\ref{contrform}. Observe that
$\Upsilon ({\cal U}({\frak g})_{\mu})=0$ for all $\mu\not=0$. In particular,
$\Upsilon(M)=\Upsilon(M_0)$ for any $\ad\frak g$-submodule $M$  
of ${\cal U}({\frak g})$.
\subsubsection{}
\begin{lem}{projann}
Let ${M}$ be an $\ad\frak g$-submodule of ${\cal U}({\frak g})$, and 
$\widetilde
{V}(\lambda)$ a simple $\frak g$-module of highest weight 
$\lambda$, then ${M}\subset \Ann_{{\cal U}({\frak g})}
\widetilde{V}(\lambda)\Longleftrightarrow
\Upsilon({M_0})(\lambda)=0$.
\end{lem}

\begin{pf} 
The proof is the same as that of Lemma 7.2 in~\cite{jfr}.
However, we present it here since this is a very important (and quite easy) 
fact.
Let $v_{\lambda}$ be a highest weight vector of 
$\widetilde{V}(\lambda)$ and 
${\widetilde{V}(\lambda)}_-\!:={\frak n}^-{\widetilde{V}(\lambda)}$. One has
$$\Upsilon(M)(\lambda)=0\Longleftrightarrow Mv_{\lambda}\subset 
{\widetilde{V}(\lambda)}_- .$$ 
On the other hand, from the $\ad{\frak g}$-invariance of $M$ one has
${\cal U}({\frak n}^-)M=M{\cal U}({\frak n}^-)$ and so the assumption 
$\Upsilon(M_0)(\lambda)=0$ implies
$$M\widetilde{V}(\lambda)=M{\cal U}({\frak n}^-)v_{\lambda}
={\cal U}({\frak n}^-)Mv_{\lambda}\subset {\cal U}({\frak n}^-)
{\widetilde{V}(\lambda)}_-\subset 
{\widetilde{V}(\lambda)}_-\ssubset{\widetilde{V}(\lambda)}.$$
The $\ad{\frak g}$-invariance of $M$ implies that 
${\cal U}({\frak g})M=M{\cal U}({\frak g})$ and consequently
$M\widetilde{V}(\lambda)$ is a submodule of $\widetilde{V}(\lambda)$.
Hence $M{\widetilde{V}(\lambda)}=0$.
\end{pf}

\subsection{The construction of the PRV determinants}
\label{prvcons} 
The separation theorem leads to the following contruction.
Fix a simple finite dimensional module 
$\widetilde{V}(\lambda)$ of the highest weight $\lambda\in P^+(\pi)$ 
and let $C(\lambda)$ be the corresponding isotypical component in ${\cal H}$.
\begin{equation}
\label{C}
C(\lambda)=\displaystyle\mathop{\oplus}_{j=1}^{\dim {V(\lambda)_0}} V^{(j)},
\ \ \widetilde{V}(\lambda)\iso V^{(j)}.
\end{equation}
Choose a basis $\{v_i\}$ of $\widetilde{V}(\lambda)_0$ and let  $\{v_i^j\}$
be the corresponding basis of $V^{(j)}_0$. Consider the matrix  
$PRV^{\lambda}:=(\Upsilon(v_i^j))$.
\Lem{projann} implies that the $j^{th}$ column of the matrix 
$PRV^{\lambda}$ is zero at a point $\mu\in\fh^*$
iff  $V^{(j)}\subset \Ann \widetilde{V}(\mu)$. Notice that the 
matrix $PRV^{\lambda}$ depends on the choice of the decomposition 
of the isotypical component $C(\lambda)$, on the choice of the 
basis of $\widetilde{V}(\lambda)_0$ 
and  on the choice of the isomorphisms from each $V^{(j)}$ to 
$\widetilde{V}(\lambda)$. 
For any change of these parameters, the new matrix one obtains 
is of the form $M(PRV^{\lambda})N$ where $M,N$ are invertible square 
complex matrices. Therefore, the corank of the 
$PRV^{\lambda}$ is correctly defined and
$\det PRV^{\lambda}\in {\cal S}({\frak h})$ 
is defined up to a non-zero scalar. Hence, 
for any $\mu\in {\frak h}^*$,  one has the equivalence:
$$\det {PRV}^{\lambda}(\mu)=0\Longleftrightarrow 
C(\lambda)\cap \Ann_{\cal H}\widetilde{V}(\mu)\not=0.$$ 
Moreover
\subsubsection{}
\begin{lem}{clprv} 
Take $\lambda\in  P^+(\pi)$ and $\mu\in {\frak h}^*$. Then 
\begin{enumerate}
\item $\corank PRV^{\lambda}(\mu)=[\Ann_{\cal H}\widetilde{V}(\mu):\widetilde
{V}(\lambda)]$
\item $\mu$ is a zero of 
$\det PRV^{\lambda}$ of order $\geq [\Ann_{\cal H}\widetilde{V}(\mu):
\widetilde{V}(\lambda)]$
\item for any $\lambda\in P^+(\pi)$,
the polynomial $\det PRV^{\lambda}$ is not identically zero.
\end{enumerate}
\end{lem}

\begin{pf} From the comments above (i) and (ii) are straightforward. 
The proof of (iii) coincides with the proof in 9.11.4,~\cite{jfr}. 
\end{pf}

\subsection{Case $l=1$}
\label{l1}
This case was completely described by Pinczon in~\cite{pin}.
In this case $\pi=\Delta_1^+=\{\beta\},\  
\Delta_0^+=\{2\beta\},\ \rho={1\over 2}\beta,\ {\frak h}^*={\Bbb C}\beta,
{\frak h}={\Bbb C}\varphi(\beta),\ P^+(\pi)={\Bbb N}\beta$. 

\subsubsection{}
\label{odd}
Denote by ${\widetilde V}(n)$ a simple module of the highest weight $n\beta$.
Call a simple module $V$ {\em odd} if $V\cong {\widetilde V}(n)$
for some odd number $n$.

Since $\dim {\widetilde V}(m)_0=1$ for all $m\in {\Bbb N}$ the separation 
theorem implies that 
$${\cal H}={\underset{m=0}{\overset{\infty}\oplus}}{\widetilde V}(m).$$
Taking into account that for any $m\in {\Bbb N},\ 
s_{\beta}.(m\beta)\!=\!-(m+1)\beta\not\in P^+(\pi)$, 
we conclude from~\Rem{hesrem} that
$$[H^n:\widetilde {V}(m)]=P_n(m\beta)$$
where the coefficients $P_n(m\beta)$ are given by the formula
$$\displaystyle\sum_{n=0}^{\infty}\displaystyle\sum_{m=0}^{\infty}P_n(m\beta)
e^{m\beta}q^n=
(1-qe^{2\beta})^{-1}(1+q^2e^{\beta})=1+qe^{2\beta}+
\displaystyle\sum_{n=2}^{\infty}q^n(e^{2n\beta}+e^{(2n-3)\beta}).$$
Hence 
\begin{equation}
H^n=\widetilde {V}(2n)\oplus\widetilde {V}(2n-3)\mbox{ for }n\geq 2
\end{equation}
Since ${\cal S}({\frak g})={ H}\otimes {{\cal S}({\frak g})}^{\frak g}$, 
we conclude that for $m=2n, 2n-3$ one has 
$[{\cal S}^k({\frak g}):\widetilde {V}(m)]=0$ for $0\leq k<n$
and $[{\cal S}^n({\frak g}):\widetilde {V}(m)]=1$. Hence
${\cal F}^n$ contains a unique copy of $\widetilde {V}(2n)$ and 
$\widetilde {V}(2n-3)$ which lie in ${\cal H}$. 
This shows in particular that, for any $n\geq 0$, the degrees of 
$\det PRV^{2n\beta}$ and $\det PRV^{2n+1}$ are less than or equal to  
$n$ and $n+2$ respectively.

\subsubsection{}
\label{prvl1}
Let $V$ be the copy of $\widetilde{V}(2n)$  in ${\cal H}$ and $v$ its 
highest weight vector. For any $k\geq 0$,
$V\subset \Ann_{\cal H} \widetilde{V}(k)$ iff $v\widetilde{V}(k)=0$. As
$\Omega({\widetilde{V}(k)})=\{k,k-1,\ldots,0,\ldots,-k+1,-k\}$, 
the condition $n>k$ is 
sufficient to ensure that $v\widetilde{V}(k)=0$ and therefore that 
$V\subset \Ann_{\cal H} \widetilde{V}(k)$.
So $\det PRV^{2n\beta}$ is divisible by 
$\mathop{\prod}_{k=0}^{n-1}(\varphi(\beta)-k)$.
Analogously we obtain also the divisibility of  $\det PRV^{(2n+1)\beta}$ by
$\mathop{\prod}_{k=0}^{n}(\varphi(\beta)-k)$. 
Taking into account that the degree of $\det PRV^{2n\beta}\leq n$ 
we conclude that, up to a non-zero scalar, 
$$\det PRV^{2n\beta}=\mathop{\prod}_{k=0}^{n-1}(\varphi(\beta)-k).$$ 

\subsubsection{}
The factorization of $\det PRV^{(2n+1)\beta}$ requires some extra work. 
The Casimir operator $C_0$ of ${\cal U}({\frak g}_0)$
acts on a ${\frak g}_0$-Verma module $M(\mu)$ of highest weight $\mu$ 
by the scalar $(\mu,\mu+2\rho_0)$.
Since ${\widetilde M}(\mu)=M(\mu)\oplus M(\mu-\beta)$ as a 
${\frak g}_0$-module,
it follows that $C_0$ acts on ${\widetilde M}(\mu)$ by a scalar
iff $\mu=-\beta/2=-\rho$. 
Hence $(C_0+{3\over 4})\in\Ann {\widetilde M}(-\rho)$.
Therefore $\Ann {\widetilde M}(-\rho)$ contains $V\!:=\ad {\cal U}({\frak g})
(C_0+{3\over 4})$. 
Since ${\frak g}_0$ acts trivially on $(C_0+{3\over 4})$ and 
$\Omega({\frak g}_1)=\{\beta,-\beta\}$, it follows that 
$\ V\subseteq V^{(1)}\oplus V^{(2)}$ where $V^{(1)}\cong{\widetilde V}(1),\ 
V^{(2)}\cong {\widetilde V}(0).$
On the other hand, $C_0\not\in {\cal Z}({\frak g})$ so 
$V^{(2)}\not=V$. Hence $V^{(1)}\subseteq V$.
Since  $C_0\in {\cal F}^2$, we conclude
that $\Ann {\widetilde M}(-\rho)$ contains a copy $V^{(1)}$ of 
${\widetilde V}(1)$
and this copy lies in ${\cal F}^2$. Thus $V^{(1)}$ lies in ${\cal H}$.
Fix a primitive vector $a$ of $V^{(1)}$. Taking into account that 
$X_{2\beta}$ is a primitive vector of ${\frak g}\cong{\widetilde V}(2)$  
and  ${\cal U}({\frak g})$ is a domain, we conclude that the product
$X_{2\beta}^ka$ is also a primitive vector of the weight $(2k+1)\beta$ and so
$\ad {\cal U}({\frak g})(X_{2\beta}^ka)\cong{\widetilde V}(2k+1)$.
Since $X_{2\beta}^ka\in {\cal F}^{k+2}$ and $a\in \Ann {\widetilde M}(-\rho)$, 
it follows that $\ \ad {\cal U}({\frak g})(X_{2\beta}^ka)\subset 
\Ann_{\cal H} {\widetilde M}(-\rho)$. Hence 
\begin{equation}
\label{annl1}
\Ann_{\cal H} {\widetilde M}(-\rho)\supseteq 
{\underset{k=0}{\overset{\infty}\oplus}}
{\widetilde V}(2k+1).
\end{equation}

Since $W.(-\rho)=\{-\rho\}$, the module ${\widetilde M}(-\rho)$ is simple.
So (\ref{annl1}) means that $\varphi(\beta)+{1\over 2}$ divides all
$\det PRV^{(2n+1)\beta}$, $n\geq 0$. 
We conclude from~\ref{prvl1} that, up to a non-zero scalar, 
$$\det PRV^{(2n+1)\beta}=(\varphi(\beta)+{1\over 2})
\mathop{\prod}_{k=0}^{n}(\varphi(\beta)-k).$$

From  the factorization  of $\det PRV^{2n\beta}$ and $\det PRV^{(2n+1)\beta}$ 
we derive that 
(\ref{annl1}) is an equality and that $\Ann_{\cal H} {\widetilde M}(\mu)=0$
if $\widetilde{M}(\mu)$ simple and $\mu\not=-\rho$.  
Since $W.(-\rho)=\{-\rho\}$,
$\Hom_{\frak g}(\widetilde{M}(-\rho),\widetilde{M}(\mu))=0$ for any 
$\mu\not=-\rho$. 
Hence $\Ann_{\cal H} {\widetilde M}(\mu)=0$ for any $\mu\not=-\rho$. 
We summarize: 
\begin{eqnarray*} \mu\not=-\rho & \Longrightarrow & 
\Ann_{\cal H}\widetilde M(\mu)=0\\
\mu=-\rho & \Longrightarrow & \Ann_{\cal H}\widetilde M(-\rho)
={\underset{i=0}{\overset{\infty}\oplus}}
{\widetilde V}(2i+1)
\end{eqnarray*}

In the sequel we will need the following
\subsubsection{}
\begin{lem}{lemcasel1}
Take ${\frak g}=\osp(1,2)$. For any odd submodule $V$ of 
${\cal U}({\frak g})$ (see~\ref{odd}) and any $v\in V$, the Harish-Chandra
projection $\Upsilon (v)$ is divisible by $(\varphi(\beta)+{1\over 2})$.
\end{lem}
\begin{pf}
Since ${\cal Z}({\frak g})$ acts on ${\widetilde M}(-\rho)$ by scalars,
the separation theorem implies that as ${\frak g}$-module 
$${\cal U}({\frak g})/(\Ann {\widetilde M}(-\rho))\cong{\cal H}/
(\Ann_{\cal H} {\widetilde M}(-\rho))={\underset{i=0}{\overset{\infty}\oplus}}
{\widetilde V}(2i).$$
Thus $\Ann {\widetilde M}(-\rho)$ contains all odd submodules of 
${\cal U}({\frak g})$.
By~\Lem{projann}, $\Upsilon (v)(-\rho)=0$ for all $v\in V$. 
Since $\Upsilon (v)$ is just a polynomial of one variable, it follows that
$\Upsilon (v)$ is divisible by
$\varphi(\beta)+(\rho,\beta)=\varphi(\beta)+{1\over 2}$ as required.
\end{pf} 

\subsection{A bound for the total degrees of PRV determinants}
\label{estimexotic}
Retain notations of~\ref{prvcons}. 
Fix a finite dimensional module ${\widetilde V}(\lambda)$ and
a decomposition (\ref{C}) of
its isotypical component $C(\lambda)$ in ${\cal H}$.
For each copy $V^{(i)}\subset {\cal H}$ let $n_i$ be such that
$\gr_{\cal F} V^{(i)}\subset H^{n_i}$.
The subspaces ${\cal F}^{(k)}$ are stable under the projection
$\Upsilon$ so 
$\Upsilon(V^{(i)})\subset {\cal F}^{(n_i)}\cap {\cal U}({\frak h})$.
Thus any entry $\Upsilon(v^i_j)$ of $i$-th column of $PRV^{\lambda}$
has degree (as polynomial in $S({\frak h})$) less than or equal to $n_i$.
Using prime to denote the derivative with respect to $q$ we obtain
$$\text {degree }\det PRV^{\lambda}\leq\displaystyle
\sum_{r=0}^{\infty} r[H^r:{\widetilde V}(\lambda)]=P'_{\lambda}(1)$$
where $\det PRV^{\lambda}$ is considered as a polynomial in $S({\frak h})$.

Using notations of~\ref{nthess} and~\Thm{thmhess}
we can rewrite~(\ref{eqthmhess}) as $P_{\lambda}(q)=\pi_0(Q_{\lambda}(q))$.
Hence
\begin{equation}
\label{estimth}
\text {degree }\det PRV^{\lambda}\leq\pi_0(Q_{\lambda}'(1)).
\end{equation}
Taking the derivative of~(\ref{eqthmhess})  on $q$, we obtain
$$Q_{\lambda}'(q)=Q_{\lambda}(q)\left(\displaystyle\sum_{\alpha\in\Delta_0^+}
e^{-\alpha}(1-qe^{-\alpha})^{-1}+2lq^{2l-1}e^{-\beta_1}
(1+q^{2l}e^{-\beta_1})^{-1}
\right).$$
By~\Rem{remq1}, $Q_{\lambda}(1)=\ch {\widetilde V}(\lambda)$. Thus
\begin{equation}
\label{estimprv1}
\begin{array}{lcl}
Q_{\lambda}'(1)&=&Q_{\lambda}(1)\left(\displaystyle\sum_{\alpha\in\Delta_0^+}
e^{-\alpha}(1-e^{-\alpha})^{-1}+2le^{-\beta_1}(1+e^{-\beta_1})^{-1}\right)\\
&=&\displaystyle\sum_{\alpha\in\Delta^+_0}\displaystyle
\sum_{m=1}^{\infty}e^{-m\alpha}\ch {\widetilde V}(\lambda)+
2l\displaystyle\sum_{m=1}^{\infty} (-1)^{m+1}e^{-m\beta_1}
\ch {\widetilde V}(\lambda)
\end{array}
\end{equation}
Since $\pi_0(e^{-\mu}\ch {\widetilde V}(\lambda))
=\dim {\widetilde V}(\lambda)_{\mu}\,$ for any $\mu\in P(\pi)$,
the formulas~(\ref{estimth}),~(\ref{estimprv1}) imply
\begin{equation}
\label{estimprv}
\text {degree }\det PRV^{\lambda}\leq \displaystyle\sum_{m=1}^{\infty}
\left(\displaystyle\sum_{\alpha\in\Delta^+_0} 
\dim \widetilde {V}(\lambda)_{m\alpha}+
2l(-1)^{m+1}\dim \widetilde {V}(\lambda)_{m\beta_1}\right)
\end{equation}

\subsection{Factorization of PRV determinants} Let us formulate the main 
theorem of this Section:

\begin{thm}{thmprv} For all $\lambda\in P^+(\pi)$, one has up to a non-zero 
scalar, 

$\vtop{\hbox to \textwidth{$\det
{PRV}^{\lambda}=\underbrace{\displaystyle\mathop{\prod}_
{\alpha\in\overline{\Delta_0^+}}
\mathop{\prod}_{m=1}^{\infty}
{\bigl[\varphi(\alpha)+(\rho,\alpha)-{1\over
2}m(\alpha,\alpha)\bigr]}^{\dim{\widetilde{V}(\lambda)}_{m\alpha}}}_{\mbox
{\footnotesize standard
factors}}\times$\hfill}
\hbox to\textwidth{\hfill$\underbrace{\displaystyle\mathop{\prod}
_{\alpha\in{\Delta_1^+}}
\mathop{\prod}_{m=1}^{\infty}
{\bigl[\varphi(\alpha)+(\rho,\alpha)-{1\over2}(2m-1)(\alpha,\alpha)\bigr]}^
{dim{\widetilde{V}(\lambda)}_{(2m-1)\alpha}}}_{\mbox{\footnotesize 
 standard factors}}\ \times $\hfill}
\hbox to \textwidth{\hfill$\underbrace{\displaystyle\mathop{\prod}_{\alpha\in
\Delta_1^+}
{\bigl[\varphi(\alpha)+(\rho,\alpha)\bigr]}^{\mathop{\sum}_{i=1}^{\infty}
{(-1)}^{i+1}
\dim {\widetilde{V}(\lambda)}_{i\alpha}}}_
{\mbox{\footnotesize exotic factors}}$}}$
\end{thm}

The rest of this section is devoted to the proof of this theorem.

\subsection{Standard factors of the PRV determinants}
\label{stfactor}
In the subsections~\ref{stfactor}---~\ref{pfthmLL} we shall explain 
briefly how
certain factors of the PRV determinants can be obtained
following the procedure used by Joseph in the classical case~\cite{jfr}, 
Chapter 9. These factors will be called "standard factors". 

Let first fix some notations. Let $M=M_0\oplus M_1$ be a ${\frak g}$-module.
We endow the
space $\Hom_{{\Bbb C}}(M,M)$ of ${\Bbb C}$-linear endomorphisms of $M$ 
with the natural supervector space structure. 
The ${\frak g}$-module structure on 
$\Hom_{{\Bbb C}}(M,M)$ is defined on homogenous elements $x\in {\frak g},\ 
f\in \Hom_{{\Bbb C}}(M,M)$ by the formula
$$ (xf)(m)=x(f(m))-{(-1)}^{|x||f|}f(xm)$$
The locally finite part of $\Hom_{{\Bbb C}}(M,M)$ with respect to this 
${\frak g}$-module structure
will be denoted $\F(M,M)$. For the above ${\frak g}$-module structure on
 $\Hom_{{\Bbb C}}(M,M)$ the natural map
 ${\cal U}({\frak g})\longrightarrow \Hom_{{\Bbb C}}(M,M)$ is 
a ${\frak g}$-module map for the adjoint action on ${\cal U}({\frak g})$. 
The image is contained in $\F(M,M)$.

Take $\mu\in {\frak h}^*$ and consider the induced morphism 
$${\cal H}\slash \Ann_{\cal H}\widetilde{V}(\mu)\longrightarrow 
\F(\widetilde{V}(\mu),\widetilde{V}(\mu)).$$
The idea to find the "standard zeroes" of $\det PRV^{\lambda}$ is to 
construct some $\mu$ such that the
multiplicity of $\widetilde{V}(\lambda)$ in 
$\F(\widetilde{V}(\mu),\widetilde{V}(\mu))$ is less than the multplicity of 
$\widetilde{V}(\lambda)$
in ${\cal H}$. 

Define for $\alpha\in \Delta_{irr}$,
$$\Lambda_{m,\alpha}\!:=\{\mu\in {\frak h}^*|\ \langle\mu+\rho,\alpha\rangle
=m,\ \  \langle\mu+\rho,\beta\rangle\not\in {\Bbb Z}\ ,\forall 
\beta\in \Delta^+_{irr}\backslash 
\{\alpha\}\}$$
for all $m\in {\Bbb N}^+$ if $\alpha\in \overline{\Delta_0^+}$ and for all odd 
$m$ if $\alpha\in \Delta_1^+$. For such choices of $m$ and $\alpha$, the set 
 $\Lambda_{m,\alpha}$ is obviously non-empty. 

 One has the
\subsubsection{}
\begin{thm}{thmLambda} Let $\lambda\in P^+(\pi)$, $\alpha\in \Delta^+_{irr}$, 
$m\in {\Bbb N}^+$ (assumed
to be odd if $\alpha$ is odd), $\mu\in \Lambda_{m,\alpha}$ then 
$$\bigl[\F(\widetilde{V}(\mu),\widetilde{V}(\mu)):\widetilde{V}(\lambda)\bigr]=
\dim {\widetilde{V}(\lambda)}_0- \dim {\widetilde{V}(\lambda)}_{m\alpha}.$$
\end{thm}
 
Assume $m,\alpha$ chosen as above. Then by~\Lem{clprv} and the theorem,
$\mu\in \Lambda_{m,\alpha}$ is a zero of $\det PRV^{\lambda}$ of order 
$\geq \dim {\widetilde{V}(\lambda)}_{m\alpha}$. The same holds for 
the Zariski closure of 
$\Lambda_{m,\alpha}$, namely for the hyperplane $(\mu+\rho,\alpha)-
m\displaystyle {(\alpha,\alpha)\over 2}$.
This means that $\bigl[\varphi(\alpha)+(\rho,\alpha)-m\displaystyle
{(\alpha,\alpha)\over 2}\bigr]$ divides $\det PRV^{\lambda}$. 
Consequently we obtain the 
\subsubsection{}
\begin{cor}{stfact} For any $\lambda\in P^+(\pi)$, the determinant 
$\det PRV^{\lambda}$ is divisible by 

$\vtop{\hbox to
\textwidth{$\quad \displaystyle\mathop{\prod}_{\alpha\in\overline{\Delta_0^+}}
\mathop{\prod}_{m=1}^{\infty}
{\bigl[\varphi(\alpha)+(\rho,\alpha)-{1\over
2}m(\alpha,\alpha)\bigr]}^{\dim{\widetilde{V}(\lambda)}_{m\alpha}}\times$
\hfill}
\hbox to \textwidth{\hfill$\displaystyle\mathop{\prod}_{\alpha\in{\Delta_1^+}}
\mathop{\prod}_{m=1}^{\infty}
{\bigl[\varphi(\alpha)+(\rho,\alpha)-{1\over2}(2m-1)(\alpha,\alpha)\bigr]}^
{dim{\widetilde{V}(\lambda)}_{(2m-1)\alpha}}\hfill$}}$
\end{cor}

\subsection{Sketch of the proof of~\Thm{thmLambda}}
\label{pfthmLL}  This subsection being very close to the 
proposition 9.7,~\cite{jfr}, we present only  a
sketch of the proof. First comes a    

\subsubsection{}
\begin{defn}{defdom} We call $\lambda\in {\frak h}^*$ dominant if 
$$\begin{array}{lll} \forall \alpha\in \overline{\Delta_0^+},&&
\langle\lambda+\rho,\alpha\rangle\not\in -{\Bbb N}^+ \\
\forall \alpha\in \Delta_1^+,&& \langle\lambda+\rho,\alpha\rangle\not\in 
-2{\Bbb N}-1\end{array}$$
\end{defn}
One has  the key 

\subsubsection{}
\begin{lem}{projectivity} The following conditions on $\lambda$ are
equivalent
\begin{enumerate}
\item $\widetilde{M}(\lambda)$ is projective in $\widetilde{\cal O}$
\item $\forall \mu\in{\frak h}^*,\ \Hom_{\frak g}(\widetilde{M}(\lambda),
\widetilde{M}(\mu))\not=0
\Longrightarrow \mu=\lambda$
\item $\lambda$  is dominant
\item $\forall \mu\in{\frak h}^*,\
 [\widetilde{M}(\mu):\widetilde{V}(\lambda)]\not=0
\Longrightarrow \mu=\lambda$.
\end{enumerate}
\end{lem}

One deduce from the lemma the following
 
\subsubsection{}
\begin{prop}{propdom} Let $\mu$ be dominant, $\mu'\in {\frak h}^*$ and $\lambda
\in P^+(\pi)$. Then 
$$\bigl[\F(\widetilde{M}(\mu),\widetilde{M}(\mu')):
\widetilde{V}(\lambda)\bigr]=
\dim{\widetilde{V}(\lambda)}_{\mu-\mu'}$$
\end{prop}

\subsubsection{}
Retain notations of~\Thm{thmLambda}. Recall that $\mu\in\Lambda_{m,\alpha}$ 
and hence is dominant.
According to factorization of Shapovalov determinants given 
in~\Thm{shapovalov}, $\widetilde{M}(\mu)$ contains
a unique copy of $\widetilde{M}(s_{\alpha}.\mu)$.
Moreover we have 
the short exact sequence
$$ 0\longrightarrow \widetilde{M}(s_{\alpha}.\mu)\longrightarrow 
\widetilde{M}(\mu)\longrightarrow 
\widetilde{V}(\mu)\longrightarrow 0$$ and  $\widetilde{M}(s_{\alpha}.\mu)$ 
is a simple Verma module. An easy consequence of the 
lemma~\ref{projectivity} is the  exactness of the functor
$\F(\widetilde{M}(\mu),-)$. This provides  
  the exact  sequence 
$$ 0\longrightarrow \F\bigl(\widetilde{M}(\mu),\widetilde{M}(s_{\alpha}.\mu)
\bigr)\longrightarrow 
\F\bigl(\widetilde{M}(\mu),\widetilde{M}(\mu)\bigr)\longrightarrow 
\F\bigl(\widetilde{M}(\mu),\widetilde{V}(\mu)\bigr)\longrightarrow 0$$
One the other hand, the functor $\F\bigl(-,\widetilde{V}(\mu)\bigr)$ gives
  the exact  sequence 
$$ \F\bigl(\widetilde{M}(s_{\alpha}.\mu),\widetilde{V}(\mu)\bigr)
\longleftarrow 
\F\bigl(\widetilde{M}(\mu),\widetilde{V}(\mu)\bigr)\longleftarrow 
\F\bigl(\widetilde{V}(\mu),\widetilde{V}(\mu)\bigr)\longleftarrow 0$$

An argument from Gelfand-Kirillov dimension theory shows
 that $ \F\bigl(\widetilde{M}(s_{\alpha}.\mu),\widetilde{V}(\mu)\bigr)=0$. 
The theorem then results from the proposition~\ref{propdom}.

\subsection{"Exotic" factors of PRV determinant}
Call the factors of the PRV determinants "exotic" if they are not standard.
By~\Cor{stfact}, the total degree of the product of the standard factors
of $\det PRV^{\lambda}$ is at least
$$\sum_{m=1}^{\infty}
\left(\displaystyle\sum_{\alpha\in {\overline\Delta^+}_0} 
\dim \widetilde {V}(\lambda)_{m\alpha}+
\displaystyle\sum_{\beta\in {\Delta_1^+}}
\dim \widetilde {V}(\lambda)_{(2m-1)\beta}\right).$$
Recall that ${\Delta_1^+}\subset W\beta_1$ and so 
$\dim\widetilde {V}(\lambda)_{k\beta}=\dim\widetilde {V}(\lambda)_{k\beta_1}$
for all $\beta\in {\Delta_1^+}$.
Taking into account that $\Delta^+_0=2\Delta^+_1\cup{\overline\Delta^+}_0$, 
we conclude from~(\ref{estimprv}) that
the total degree of the product of the exotic factors
of $\det PRV^{\lambda}$ is less than or equal to 
\begin{equation}
\label{estimprvex}
\sum_{\beta\in {\Delta_1^+}}\sum_{m=1}^{\infty}(-1)^{m+1}\dim 
\widetilde {V}(\lambda)_{m\beta}.
\end{equation}

\subsubsection{}\label{remprvl1}
\begin{rem}{} Consider the case $l=1$, i.e. ${\frak g}=\osp(1,2)$ 
and retain notations of~\ref{l1}. 
Then the formula (\ref{estimprvex}) takes the values 
\begin{equation}
\label{l1prvex}
\sum_{m=1}^{\infty}(-1)^{m+1}\dim {\widetilde V}(n)_{m\beta}=
\left\{
\begin{array}{ll}
0,&\ n \ \text { is even}\\
1,&\ n \ \text { is odd}
\end{array}
\right.
\end{equation}
\end{rem}

\subsection{}
\label{genl}
Let ${\frak p}$ be 
$${\frak p}\!:={\frak g}_{\pm\beta_l}\oplus {\frak g}_{\pm2\beta_l}\oplus
{\Bbb C}\varphi(\beta_l).$$
Then ${\frak p}$ is a Lie subsuperalgebra of ${\frak g}$ and 
there is an obvious isomorphism ${\frak p}\cong\osp(1,2)$ which 
maps the element $\varphi(\beta_l)$ to $\varphi(\beta)$. 

We shall denote a simple 
${\frak p}$-module of the highest weight $n\beta_l$ by ${\breve V}(n)$ and
call a simple ${\frak p}$-module $V$ "odd" if $V\cong{\breve V}(n)$ for some
odd number $n$. 
Lemma~\ref{lemcasel1} has the following useful generalization.

\subsubsection{}
\begin{lem}{lembeta}
For any odd $\ad {\frak p}$-submodule $V$ of ${\cal U}({\frak g})$ 
and any $v\in V$, the Harish-Chandra
projection $\Upsilon (v)$ is divisible by $\varphi(\beta_l)+{1\over 2}=
(\varphi(\beta_l)+(\rho,\beta_l))$.
\end{lem}
\begin{pf}
Since $\beta_l$ is a simple root of $\Delta$,  
$({\frak p}+{\frak b}^{\pm})$ are 
Lie subsuperalgebras and even parabolic subsuperalgebras of ${\frak g}$.  
Let ${\frak h}'$ be the linear span of $\{\varphi(\beta_i)\}_{i=1}^{l-1}$.
The superalgebra $({\frak p}\oplus{\frak h}')$ is the Levi factor
of the parabolic superalgebras  $({\frak p}+{\frak b}^{\pm})$.
Let ${\frak m}^{\pm}$ be the nilradical of $({\frak p}+{\frak b}^{\pm})$.
Denote by  ${\cal U}({\frak p}+{\frak h}')$ the universal enveloping 
algebra of $({\frak p}+{\frak h}')$.
Let $P$ be the projection of ${\cal U}({\frak g})$ onto 
${\cal U}({\frak p}+{\frak h}')$
with respect to decomposition 
${\cal U}({\frak g})={\cal U}({\frak p}+{\frak h}')\oplus
({\frak m}^-{\cal U}({\frak g})+{\cal U}({\frak g}){\frak m}^+).$
Then $\Upsilon=\Upsilon\circ P$. 
Fix a simple odd $\ad {\frak p}$-module $V\subset {\cal U}({\frak g})$ and 
$v\in V$. Since
$P$ is an $\ad {\frak p}$-homomorphism, it follows that either $P(V)=0$ or
$P(V)\simeq V$. In the first case the assertion obviously holds. Consider the
second case. Since
${\frak p}$ acts trivially on ${\frak h'}$ and on its universal enveloping 
algebra ${\cal U}({\frak h'})$, the multiplication map induces an isomorphism
${\cal U}({\frak p}+{\frak h}')={\cal U}({\frak p})\otimes 
{\cal U}({\frak h'})$ as ${\frak p}$-modules.
Write $P(v)=\sum p_ih_i$ with non-zero $p_i\in {\cal U}({\frak p}),\ 
h_i\in {\cal U}({\frak h'})$. Since $\ad {\cal U}({\frak p})P(v)\!\simeq\!V$ 
is a simple ${\frak p}$-module, one can suppose that 
$V\!\simeq\!\ad {\cal U}({\frak p})(p_i)$ for all $i$. 
Since $p_i\in {\cal U}({\frak p})$, $\Upsilon(p_i)$ 
is divisible by $(\varphi(\beta_l)+{1\over 2})$ for all $i$ by~\Lem{lemcasel1}.
Thus $\Upsilon(v)=\Upsilon(P(v))=\sum \Upsilon(p_ih_i)=\sum \Upsilon(p_i)h_i$ 
is divisible by $(\varphi(\beta_l)+{1\over 2})$ as required. 
This establishes the lemma.
\end{pf}

\subsubsection{}
\label{casebetan}
Fix $\lambda\in P^+(\pi)$ and an $\ad{\frak g}$-submodule 
$V\cong \widetilde {V}(\lambda)$ of ${\cal H}$. Then
there is an $\ad{\frak p}$-submodule of $V$ given by the formula
$${\breve V}\!:={\underset{m\in {\Bbb Z}}\oplus} V_{m\beta_l}.$$
Fix a decomposition of ${\breve V}$ into the sum of simple 
$\ad{\frak p}$-modules.
Taking into account that $V_0\subset {\breve V}$ and that
the dimension of the zero weight space of any simple ${\frak p}$-module
is equal to one, we conclude that the number of simple $\ad{\frak p}$-modules
in the decompositon of ${\breve V}$ is equal to $\dim V_0$. Hence
$${\breve V}={\underset{i=1}{\overset{\dim V_0}\oplus}} {\breve V}^{(i)}.$$
Using~\ref{remprvl1} we obtain
\begin{equation}
\label{alsum}
\displaystyle\sum_{m=1}^{\infty}(-1)^{m+1}
\dim {\widetilde{V}(\lambda)}_{m\beta_l}=
\displaystyle\sum_{i=1}^{\dim {V}_0} \displaystyle\sum_{m=1}^{\infty}
(-1)^{m+1}\dim {\breve V}^{(i)}_{m\beta_l}=
\# \{i:\ {\breve V}^{(i)} \text { is odd }\}.
\end{equation}

Choose a basis $\{v_i\}_{i=1}^{\dim V_0}$ of the subspace $V_0$ such
that $v_i\in {\breve V}^{(i)}$ for all $i=1,\ldots,\dim V_0$. Consider
a matrix $PRV^{\lambda}$ constructed using this basis. By~\Lem{lembeta} 
the polynomial $\Upsilon (v_i)$  is divisible by 
$(\varphi(\beta_l)+(\rho,\beta_l))$ if ${\breve V}^{(i)}$ is an odd 
${\frak p}$-module. 
Thus all entries of the $i$-th line of the  matrix $PRV^{\lambda}$ are  
divisible by $(\varphi(\beta_l)+(\rho,\beta_l))$ if ${\breve V}^{(i)}$ 
is an odd ${\frak p}$-module. Using the equality~(\ref{alsum}) we obtain the

\subsubsection{}\begin{cor}{corbeta}

(i) For any $\lambda\in P^+(\pi)$ and any $\mu\in {\frak h}^*$ such that
$(\mu+\rho,\beta_l)=0$ one has
$$\corank PRV^{\lambda}(\mu)\geq \sum_{m=1}^{\infty}(-1)^{m+1}
\dim {\widetilde V}(\lambda)_{m\beta_l}.$$

(ii) For any $\lambda\in P^+(\pi)$, the determinant $\det PRV^{\lambda}$ is
divisible by 
$$\left(\varphi(\beta_l)+(\rho,\beta_l)\right)^{\sum_{m=1}^{\infty}(-1)^{m+1}
\dim {\widetilde V}(\lambda)_{m\beta_l}}$$
\end{cor}
 
\subsection{Democracy principle}
\label{democrasy}
In~\Cor{corbeta} we found 
$m_{\lambda}\!:=\sum_{m=1}^{\infty}(-1)^{m+1}
\dim {\widetilde V}(\lambda)_{m\beta_l}$
linear exotic factors of $\det PRV^{\lambda}$. By~(\ref{estimprvex}) 
the number of linear
exotic factors is at most $l\cdot m_{\lambda}$. In order to obtain
the rest of the factors we are going to prove in this subsection 
the following democracy theorem establishing 
"equality of the odd positive roots 
according Law".

\begin{thm}{thmdem}

(i) For any $\lambda\in P^+(\pi)$ and any $\mu\in {\frak h}^*$ such that
$(\mu+\rho,\beta_i)=0$ one has
$$\corank PRV^{\lambda}(\mu)\geq \sum_{m=1}^{\infty}(-1)^{m+1}
\dim {\widetilde V}(\lambda)_{m\beta_l}.$$

(ii) For any $\lambda\in P^+(\pi)$, $\beta_i\in {\Delta_1^+}$
the determinant $\det PRV^{\lambda}$ is divisible by 
$$\left(\varphi(\beta_i)+(\rho,\beta_i)\right)^{\sum_{m=1}^{\infty}(-1)^{m+1}
\dim {\widetilde V}(\lambda)_{m\beta_l}}.$$
\end{thm} 

\subsubsection{}
Let $s_i\in W$ be the simple reflection with respect to the simple root 
$\alpha_i$,
where $\alpha_i=(\beta_i-\beta_{i+1})$ for $i<l$ and $\alpha_l=\beta_l$. 
In order to prove the theorem above we need the following lemma. 

\begin{lem}{estimN}
Fix $i\in\{1,\ldots,l-1\}$ and $n\in {\Bbb N}^+$. Take $\mu\in {\frak h}^*$
such that 
$$
\begin{array}{c}
(\beta_i,\mu+\rho)=n,\ (\beta_{i+1},\mu+\rho)=0,\\
\text {the collection}\{(\beta_k,\mu+\rho)\}_{ k\not=i+1}\ 
\text { is linearly independent over }{\Bbb Q}.
\end{array}
$$

(i) Let $N$ be the kernel of the surjective map ${\widetilde M}(\mu)\to 
{\widetilde V}(\mu)$.
Then 
$$\underset{\xi\in\Omega (N)}\min (\mu-\xi,\omega_i)=n.$$

(ii) Assume that $\lambda\in P^+(\pi)$ is such that $(\lambda,\omega_i)<n$.
Then 
$$\overline{C}(\lambda)\cap \Ann {\widetilde V}(\mu)\subseteq
\overline{C}(\lambda)\cap \Ann {\widetilde V}(s_i.\mu).$$
where $\overline{C}(\lambda)$ is the isotypical component of 
$\widetilde{V}(\lambda)$ in ${\cal U}({\frak g})$.
\end{lem} 
\begin{pf}
(i) The assumption on $\mu$ implies that 
$\langle\mu+\rho,\alpha\rangle\in {\Bbb N}^+$ only when 
$\alpha\!=\!\beta_i\pm\beta_{i+1}$ or $\alpha\!=\!\beta_{i}$. 
Since $\langle\mu+\rho,\beta_{i}\rangle=2n\not\in 2{\Bbb N}+1$
and $\langle\mu+\rho,\beta_{i}\pm\beta_{i+1}\rangle=n$, we
conclude from~\ref{weightsub} that
$$\mu-\Omega (N)=\{n(\beta_i-\beta_{i+1})+{\Bbb N}(\pi)\}\cup
\{n(\beta_i+\beta_{i+1})+{\Bbb N}(\pi)\}.$$
Hence 
$$\underset{\xi\in\Omega (N)}\min (\mu-\xi,\omega_i)=
\underset{\nu\in {\Bbb N}(\pi)}\min \{(n(\beta_i-\beta_{i+1})+\nu,\omega_i),
(n(\beta_i+\beta_{i+1})+\nu,\omega_i)\}=n$$
as required.

(ii) Since $(\mu+\rho,\beta_i-\beta_{i+1})\in {\Bbb N}^+$ it follows that
${\widetilde V}(s_i.\mu)$ is a subquotient of ${\widetilde M}(\mu)$. Therefore
$\Ann {\widetilde M}(\mu)\subseteq \Ann {\widetilde V}(s_{\alpha_i}.\mu)$.
Thus it is sufficient to show that 
$$\overline{C}(\lambda)\cap \Ann {\widetilde V}(\mu)=
\overline{C}(\lambda)\cap \Ann {\widetilde M}(\mu).$$
Assume the contrary, namely, that there exists an
$\ad {\frak g}$-submodule $V\cong {\widetilde V}(\lambda)$ of 
$\Ann {\widetilde V}(\mu)$ such that $V{\widetilde M}(\mu)\not=0$.
Since $V$ is $\ad {\frak g}$-invariant one has
$V{\cal U}({\frak g})={\cal U}({\frak g})V$. Consequently, 
$V{\widetilde M}(\mu)\not=0$ forces $Vv_{\mu}\not=0$, where
$v_{\mu}$ is a highest weight vector of ${\widetilde M}(\mu)$.
Recall from~\ref{duality} that $\Omega (V)=-\Omega (V)$.
Therefore
\begin{equation}\label{omegan}
\underset{\xi\in\Omega (Vv_{\mu})}\max (\mu-\xi,\omega_i)=
\underset{\nu\in\Omega (V)}\max (-\nu,\omega_i)=
\underset{\nu\in\Omega (V)}\max (\nu,\omega_i)<n.
\end{equation}
On the other hand, $V{\widetilde V}(\mu)=0$ and so $Vv_{\mu}\subseteq N$.
Thus~(\ref{omegan}) contradicts to (i).
\end{pf}

\subsubsection{Proof of~\Thm{thmdem}}
Fix $\lambda\in P^+(\pi)$ and set 
$$r=\sum_{m=1}^{\infty}(-1)^{m+1}\dim {\widetilde V}(\lambda)_{m\beta_l}.$$ 
For $k=1,\ldots,l$ let $x_k$ be the element of ${\cal S}({\frak h})$ given by 
the formula $x_k\!:=\varphi(\beta_k)+(\rho,\beta_k)$. 
Recall that $s_i$, $i=1,\ldots,l-1$,  is the simple reflection with respect to
the root $(\beta_i-\beta_{i+1})$ and consequently one has, for any 
$i=1,\ldots,l-1$
\begin{equation}
\label{xksi}
x_i(s_i.\mu)=x_{i+1}(\mu),\ x_{i+1}(s_i.\mu)=x_i(\mu),\ \ 
x_k(s_i.\mu)=x_k(\mu) \text { for } k\not=i,i+1.
\end{equation}
Each entry $a_{i,j}$ of the matrix $PRV^{\lambda}$ is an element
of ${\cal S}({\frak h})$ and so can be considered
as a polynomial in $x_1,\ldots , x_l$. 
We have to prove that for any $i=1,\ldots,l$ and for any $\mu$ such that 
$x_i(\mu)=0$, the corank of the matrix 
$PRV^{\lambda}(x_1(\mu),\ldots ,x_l(\mu))$ is at least $r$. 
We prove this assertion by induction on $i$.
By~\Cor{corbeta}, the assertion holds for $i=l$.
Assume that the assertion holds for some $i+1:\ 1<i<l$ and deduce it for $i$. 
Recall~\Claim{clprv}: 
\begin{equation}
\label{anncorank}
[\Ann_{{\cal H}} {\widetilde V}(\mu):{\widetilde V}(\lambda)]=
\corank PRV^{\lambda}(\mu)
\end{equation}
Thus our assumption implies that for any $\mu$ such that
$(\mu+\rho,\beta_{i+1})=0$ one has
\begin{equation}
\label{annsi}
[\Ann_{{\cal H}} {\widetilde V}(\mu):{\widetilde V}(\lambda)]\geq r.
\end{equation}
Fix, for a moment, a collection of complex numbers $\{a_k\}_{k=1}^{l}$ 
such that 
the collection $\{1,a_k:k=1,\ldots,l\}$ is linearly independent over 
${\Bbb Q}$.
For every $n\in {\Bbb N}^+$, $n>(\lambda,\omega_i)$, take 
$\mu_n\in {\frak h}^*$
such that $x_i(\mu_n)=n,\ x_{i+1}(\mu_n)=0,\ x_k(\mu_n)=a_k$ for 
$k\not=i,i+1$.
Then~\Lem{estimN} (ii) combined with the inequality~(\ref{annsi}) implies that
\begin{eqnarray*}
r\leq [\Ann_{{\cal H}} {\widetilde V}(s_i.\mu_n):{\widetilde V}(\lambda)]&=&
\corank PRV^{\lambda}(s_i.\mu_n)\\
&=&\corank PRV^{\lambda}(x_1(s_i.\mu_n),\ldots ,x_l(s_i.\mu_n)).
\end{eqnarray*}
Thus  $\corank PRV^{\lambda}$  is greater than or equal to $r$ at the points 
$$(x_1(s_i.\mu_n),\ldots,x_l(s_i.\mu_n))=
(a_1,\ldots,a_{i-1},0,n,a_{i+2},\ldots, a_l)\ \ \mbox{ (see~(\ref{xksi}))}$$
of the line $$(a_1,\ldots,a_{i-1},0,{\Bbb C},a_{i+2},\ldots, a_l).$$
Therefore $\corank PRV^{\lambda}\geq r$ on the whole line 
$(a_1,\ldots,a_{i-1},0,{\Bbb C},a_{i+2},\ldots, a_l)$.
Observe that the union of the lines 
$(a_1,\ldots,a_{i-1},0,{\Bbb C},a_{i+2},\ldots, a_l)$ such
that $\{1,a_k:k=1,\ldots,l\}$ is linearly independent over ${\Bbb Q}$ is dense
in the hyperplane $x_i=0$. Hence $\corank PRV^{\lambda}$  
is greater than or equal to 
$r$ in the whole hyperplane $x_i=0$ as required. 
This establishes~\Thm{thmdem}.

\subsection{}
\label{endpfprv}
Combining the results of
~\Cor{stfact},~\Thm{thmdem} and the bound~(\ref{estimprvex}) 
we  obtain~\Thm{thmprv}.

\section{Main theorem}
\subsection{}
\begin{thm}{main} Take $\lambda\in {\frak h}^*$. Then 
\begin{equation}
\label{maineq} 
\Ann_{{\cal U}({\frak g})}{\widetilde M}(\lambda)=
{\cal U}({\frak g})\Ann_{{\cal Z}({\frak g})}{\widetilde M}(\lambda)\  \ 
\Longleftrightarrow \ \ \forall\alpha\in \Delta_1^+:\ 
(\lambda+\rho,\alpha)\not=0.
\end{equation}
\end{thm}

The standard factors of~\Thm{thmprv} are exactly the factors which appear in 
the factorization of Shapovalov determinants---see~\Thm{shapovalov}.
Since a Verma module $\widetilde{M}(\mu)$ is simple iff $\mu$ is
not a root of any Shapovalov determinant, this proves the 
equivalence~(\ref{maineq}) for simple Verma modules. By~\ref{redtosimple},
it gives immediatly the implication "$\Leftarrow$" of~(\ref{maineq}).

\subsection{}
\label{otherim}
Let us prove the other implication.
Consider ${\frak g}_1$ as a ${\frak g_0}$-module. All its
weight-spaces are one-dimensional,
$\Omega ({\frak g}_1)=\{\pm\beta_1,\ldots,\pm\beta_l\}$ so 
${\frak g}_1\cong V(\omega_1)$. Consider the
natural representation $\widetilde{V}(\omega_1)$ of $\osp(1, 2l)$.
It is easy to see that 
$\widetilde{V}(\omega_1)=V(\omega_1)\oplus V(0)$ as ${\frak g_0}$-modules.
Again all weight-spaces of $\widetilde{V}(\omega_1)$ are one-dimensional 
and $\Omega (\widetilde{V}(\omega_1))=\{0,\pm\beta_1,\ldots,\pm\beta_l\}$.
Hence, by~\Thm{thmhess},
 ${\cal H}$ contains exactly one copy $V$ of $\widetilde{V}(\omega_1)$
which occurs in degree $2l$.
Fix $i\in\{1,\ldots,l\}$. One has
$$\mathop{\sum}_{m=1}^{\infty}{(-1)}^{m+1}
\dim{\widetilde{V}(\omega_1)}_{m\beta_i}=1$$ 
and so, by~\Thm{thmprv}, $\det PRV^{\omega_1}(\mu)=0$  for any 
$\mu\in {\frak h}^*$
such that $(\mu+\rho,\beta_i)=0$. By~\Lem{clprv}, it implies that 
$V\subset \Ann_{{\cal U}({\frak g})}\widetilde{V}(\mu)$. 
In other words, $V\subset \Ann_{{\cal U}({\frak g})}\widetilde{M}(\mu)$ 
for any $\mu$ in the set 
$$\left\{\mu\in {\frak h}^*|\ (\mu+\rho,\beta_i)=0 \mbox { and }
 \widetilde{M}(\mu)\mbox{ is simple}\right\}.$$
According to the criterion of simplicity given in~\Cor{simpleV}, this set is 
Zariski dense in the hyperplane
$(\mu+\rho,\beta_i)=0$. 

The following lemma proves that in fact 
$V\subset \Ann_{{\cal U}({\frak g})}\widetilde{M}(\mu)$ 
for any $\mu$ such that $(\mu+\rho,\beta_i)=0$. 
This gives the implication "$\Longrightarrow$" of~(\ref{maineq}).

\subsection{}
\begin{lem}{zarcl} 
Let $V$ be an $\ad$-invariant subspace of ${\cal U}({\frak g})$. Then the set
 
$\left\{\lambda\in{\frak h}^*,\ V\subset\Ann_{{\cal U}({\frak g})} 
\widetilde{M}(\lambda)\right\}$ is a Zariski closed subset of ${\frak h}^*$.
\end{lem}{}

\begin{pf} For any $\lambda\in {\frak h}^*$ denote by $v_{\lambda}$
a highest weight vector of $\widetilde{M}(\lambda)$.
Remark that the $\ad$-invariance of $V$ implies that 
${\cal U}({\frak g})V=V{\cal U}({\frak g})$.
Hence, for any $\lambda\in {\frak h}^*$, one has the equivalence $V\subset 
\Ann_{{\cal U}({\frak g})}\widetilde{M}(\lambda)
\Longleftrightarrow V.v_{\lambda}=0$. Let 
 $\{a_j\}_{j\in J}$ be a basis of $V$. Then 
$$\left\{\lambda\in{\frak h}^*,\ V\subset\Ann_{{\cal U}({\frak g})} 
\widetilde{M}(\lambda)\right\}=\mathop{\bigcap}_{j\in J}
\left\{\lambda\in {\frak h}^*,\ x.v_{\lambda}=0\right\}.$$

Thus it is sufficient to show that for every $a\in {\cal U}({\frak g})$,
the set $\left\{\lambda\in {\frak h}^*,\ a.v_{\lambda}=0\right\}$
is a Zariski closed subset of ${\frak h}^*$.
Using the triangular decomposition ${\cal U}({\frak g})=
{\cal U}({\frak n}^-)\otimes{\cal U}({\frak h})\otimes{\cal U}({\frak n}^+)$,
write
$a=y_1P_1+\ldots+y_rP_r+x$ where $y_1,\ldots,y_r$ are linearly independent
elements of ${\cal U}({\frak n}^-)$, $P_1,\ldots,P_r\in{\cal U}({\frak h})$ 
and $x\in {\cal U}({\frak g}){\frak n}^+$. Then
$$\begin{array}{lcl} a.v_{\lambda}=0 &\Longleftrightarrow & \bigl(\displaystyle
\mathop{\sum}_i y_iP_i\bigr).v_{\lambda}=0\\
&\Longleftrightarrow & \displaystyle\mathop{\sum}_i P_i(\lambda)y_i.
v_{\lambda}=0\\
&\Longleftrightarrow & \forall i\ P_i(\lambda)=0.\
\end{array}.$$
The lemma is proved. 
\end{pf}
This completes the proof of~\Thm{main}.

%%%%%%%%%%%%%%%%  biblio1.tex


\begin{thebibliography}{MMMMM}

\bibitem[AL]{al} M.~Aubry and J.-M. Lemaire, 
Zero divisors in enveloping algebras
of graded Lie algebras, J. Pure Appl. Algebra {\bf 38} (1985), p. 159---166.

\bibitem[D]{di} J.~Dixmier, Alg\`ebres enveloppantes, Cahiers scientifiques, 
vol. {\bf 37}, Gauthier-Villars, Paris (1974). 

\bibitem[He]{he} W.~H.~Hesselink, Characters of the nullcone, 
Math. Ann. {\bf 252}, (1980) p.179---182.

\bibitem[H]{h} J.~E.~Humphreys, Reflection groups and Coxeter groups, 
Cambridge studies in advanced mathematics, {\bf 29}.

\bibitem[Jak]{jak} H.~P.~Jakobsen, The full set of unitarizable
highest weight modules of basic classical Lie superalgebras, Memoirs
of the Amer. Math. Society, {\bf 532} (1994).

\bibitem[J1]{jfr} A.~Joseph, Sur l'annulateur d'un module de Verma,
NATO Adv. Sci. Inst. Ser. C Math. Phys. Sci., {\bf 514} (1997), p. 237--300.


\bibitem[J2]{jq} A.~Joseph, Quantum groups and their primitive ideals, 
Springer Verlag  (1995).

\bibitem[JL]{jl} A.~Joseph  and G.~Letzter, Verma modules
annihilators and quantized enveloping algebras, Ann. Ec. Norm. Sup.
s\'erie 4, t. {\bf 28} (1995), p.493---526.

\bibitem[K1]{k1} V.~G.~Kac, Characters of typical 
representations of classical Lie superalgebras, 
Comm. in Algebra, vol. {\bf 5}, (1997), p.889---897.

\bibitem[K2]{k2} V.~G.~Kac, Lie Superalgebras, Advances in Math.
{\bf 26}, (1977), p.8---96.

\bibitem[K3]{k3} V.~G.~Kac, 
 Highest weight representations of conformal current algebras,
   Topological and geometrical methods in field theory, 
   Proc. Symp., Espoo/Finl. 1986, {\bf 3-15} (1987).


\bibitem[Mu1]{mu} Ian.~M.~Musson, On the center of the enveloping
algebra of a classical simple Lie superalgebra, J. of Algebra,
{\bf 193} (1997), p.75---101.

\bibitem[Mu2]{mu2} Ian.~M.~Musson, The enveloping algebra of the Lie
superalgebra $\osp(1,2r)$, Rep. Theory, Vol {\bf 1} (1997), 405-423.

\bibitem[Pi]{pin} G.~Pinczon, The enveloping algebra of the Lie
superalgebra $\osp (1,2)$, J. of Algebra,
{\bf 132} (1990), p.219---242.

\bibitem[Sch]{sch} M.~Scheunert, The theory of Lie superalgebras,
Lect. Notes in Math., {\bf 716}, Springer-Verlag (1979).





\end{thebibliography}
\end{document}